\documentclass[a4paper,10pt]{article}
\usepackage{amsmath,amsthm,latexsym,makeidx,amscd,amssymb}
\usepackage[all]{xy}

\usepackage{graphics}

\makeindex

 \numberwithin{equation}{section}

\newtheorem{prop}{Proposition}[section]
\newtheorem{lem}[prop]{Lemma}
\newtheorem{ddd}[prop]{Definition}
\newtheorem{theorem}[prop]{Theorem}
\newtheorem{cor}[prop]{Corollary}

\newcommand{\DF}{\mathop{\mbox{\rm DF}}} 

\newcommand{\End}{{\rm End}}

\newcommand{\ind}{\mbox{\rm ind}}
\newcommand{\Ind}{\mathop{\mbox{\rm Ind}}}

\newcommand{\dom}{\mathop{\rm dom}}

\newcommand{\U}{{\cal U}}

\newcommand{\A}{{\cal A}}
\newcommand{\B}{{\cal B}}

\newcommand{\C}{C^{\infty}}

\newcommand{\ra}{\partial}

\newcommand{\ten}{\otimes}

\newcommand{\ve}{\varepsilon}
\newcommand{\ov}{\overline}

\newcommand{\dira}{\partial\!\!\!/}

\DeclareMathOperator{\spfl}{sf}
\DeclareMathOperator{\Spfl}{Sf}
\DeclareMathOperator{\supp}{supp}
\DeclareMathOperator{\Ran}{Ran}
\DeclareMathOperator{\Ker}{Ker}

\DeclareMathOperator{\diag}{diag}
\DeclareMathOperator{\Mod}{mod}
\DeclareMathOperator{\ev}{ev}
\DeclareMathOperator{\Ima}{Im}

\def\bbbr{{\rm I\!R}} 
\def\bbbn{{\rm I\!N}} 

\def\bbbc{{\rm I\!C}}

\def\bbbq{{\mathchoice {\setbox0=\hbox{$\displaystyle\rm Q$}\hbox{\raise
0.15\ht0\hbox to0pt{\kern0.4\wd0\vrule height0.8\ht0\hss}\box0}}
{\setbox0=\hbox{$\textstyle\rm Q$}\hbox{\raise
0.15\ht0\hbox to0pt{\kern0.4\wd0\vrule height0.8\ht0\hss}\box0}}
{\setbox0=\hbox{$\scriptstyle\rm Q$}\hbox{\raise
0.15\ht0\hbox to0pt{\kern0.4\wd0\vrule height0.7\ht0\hss}\box0}}
{\setbox0=\hbox{$\scriptscriptstyle\rm Q$}\hbox{\raise
0.15\ht0\hbox to0pt{\kern0.4\wd0\vrule height0.7\ht0\hss}\box0}}}}

\def\bbbz{{\mathchoice {\hbox{$\sf\textstyle Z\kern-0.4em Z$}}
{\hbox{$\sf\textstyle Z\kern-0.4em Z$}}
{\hbox{$\sf\scriptstyle Z\kern-0.3em Z$}}
{\hbox{$\sf\scriptscriptstyle Z\kern-0.2em Z$}}}}

\def\bbbc{{\mathchoice {\setbox0=\hbox{$\displaystyle\rm C$}\hbox{\hbox
to0pt{\kern0.4\wd0\vrule height0.9\ht0\hss}\box0}}
{\setbox0=\hbox{$\textstyle\rm C$}\hbox{\hbox
to0pt{\kern0.4\wd0\vrule height0.9\ht0\hss}\box0}}
{\setbox0=\hbox{$\scriptstyle\rm C$}\hbox{\hbox
to0pt{\kern0.4\wd0\vrule height0.9\ht0\hss}\box0}}
{\setbox0=\hbox{$\scriptscriptstyle\rm C$}\hbox{\hbox
to0pt{\kern0.4\wd0\vrule height0.9\ht0\hss}\box0}}}}

\makeindex
\title{On the noncommutative spectral flow}
\author{Charlotte Wahl\footnote{This research was supported by a grant of AdvanceVT}}
\date{}

\begin{document}
\maketitle

\begin{abstract}
We define and study the noncommutative spectral flow for paths of regular selfadjoint Fredholm operators on a Hilbert $C^*$-module. We give an axiomatic description and discuss some applications. One of them is the definition of a noncommutative Maslov index for paths of Lagrangians, which appears in a splitting formula for the spectral flow. Analogously we study the spectral flow for odd operators on a $\bbbz/2$-graded module. 
\end{abstract}

\section{Introduction}

The spectral flow as introduced in \cite{aps} assigns to a continuous path of selfadjoint Fredholm operators on a Hilbert space the net number of eigenvalues changing sign along the path. It is homotopy invariant and closely related to Bott periodicity. Motivated by geometric applications it has been generalized in several directions. The purpose of this paper is to unify and generalize two of them -- the spectral flow for unbounded selfadjoint Fredholm operators and the noncommutative spectral flow, which has been defined for paths of selfadjoint operators with compact resolvents, -- by defining a spectral flow for paths of regular selfadjoint Fredholm operators on the standard Hilbert $\A$-module $H_{\A}$, where $\A$ is a unital $C^*$-algebra. In particular this includes paths of bounded selfadjoint operators on $H_{\A}$.

The most general definition of the spectral flow for selfadjoint Fredholm operators on a separable Hilbert space so far has been given by Booss-Bavnek--Lesch--Phillips \cite{blp}. Here the path is assumed to be continuous in the gap metric, i.e. the resolvents are assumed to depend continuously on the parameter. The continuity condition we impose will be weaker than gap continuity. For several reasons this is desirable: First of all the condition is easier to verify. Also there are elementary examples of paths of elliptic operators on a noncompact manifold with coefficients depending continuously on the parameter which are not gap continuous but for which a spectral flow can be defined. From a conceptual viewpoint the spectral flow should be invariant under conjugation by a strongly continuous path of unitaries such that it is well-defined also if the Hilbert space depends on the parameter, and it should be invariant under taking the bounded transform $D \mapsto D(1+D^2)^{-1/2}$. \\
See \cite{le} for an axiomatic definition and an overview over the development of the classical spectral flow. 

Dai--Zhang introduced a spectral flow for paths of families of elliptic pseudodifferential operators of order one, which has its values in $K^0$ of the base space \cite{dz2}. A modification is the noncommutative spectral flow defined for example for Dirac operators associated to $C^*$-vector bundles \cite{wu}\cite{lp2}. Both definitions rely on the concept of spectral sections developed by Melrose--Piazza \cite{mp}. Wu generalized the concept of spectral sections  to selfadjoint operators with compact resolvents on a Hilbert $C^*$-module \cite{wu}.  The noncommutative spectral flow has been applied to the study of higher signatures and higher transverse signatures of foliated bundles \cite{lp2} \cite{lp3}.

Motivating examples for our generalization are paths of elliptic operators on a noncompact manifold as well as paths of well-posed elliptic boundary value problems (see \cite{blp}) and the corresponding problems for families of operators resp. operators over $C^*$-algebras. We also expect applications to geometric operators on noncommutative spaces (for example as in \cite{pr}). Another motivation comes from the pairing of odd $K$-theory with $KK$-theory, which we expressed in terms of the noncommutative spectral flow in \cite{wa2} using results of the present paper.
 
Two classical approaches are the basis for our definition of the spectral flow: Loosely speaking, the first defines the spectral flow as a measure of the change of appropriate (generalized) spectral projections along the path. The second associates a loop of unitaries to the path to which a winding number is applied. Definitions of the spectral flow using the first approach are typically easy to identify with the intuitive meaning of the spectral flow. In the noncommutative context the second approach is more natural. Here Bott periodicity provides a generalization of the winding number. The family spectral flow and the noncommutative spectral flow have been defined via the first approach with spectral sections as a generalization of spectral projections. The relation to Bott periodicity has been studied in Wu's unpublished preprint \cite{wu}. For the definition of the spectral flow for unbounded selfadjoint Fredholm operators on a Hilbert space both methods have been applied \cite{blp}.  

The concept of spectral sections uses the compactness of the resolvents in a crucial way. One purpose of this paper is to develop an appropriate generalization for operators with noncompact resolvents. For that aim we introduce a relative index for pairs of projections on $H_{\A}$ differing by a compact operator, which generalizes the difference element of spectral sections. This notion is well-understood on a Hilbert space but has not been studied before in a systematic way on a Hilbert $C^*$-module. A priori there are several natural candidates for a relative index. One of our results is that the relative index is uniquely defined by a few axioms. Applied to families the relative index is strongly related to the generalized Maslov index introduced in \cite{n}. However our approach is different and seems to be more appropriate in the context of Hilbert $C^*$-modules.

The concept of generalized spectral sections is then used to give a definition of the noncommutative spectral flow in analogy to the family spectral flow in \cite{dz2}. However, generalized spectral sections need not exist for a given path. The general definition of the noncommutative spectral flow is in terms of a unitary transform and Bott periodicity. It applies to any path of regular selfadjoint Fredholm operators $(D_t)_{t \in [0,1]}$ with invertible $D_0$ and $D_1$ such that $(D_t)_{t \in [0,1]}$ is regular and Fredholm as an operator on the Hilbert $C([0,1],\A)$-module $C([0,1],H_{\A})$. It works also for $\sigma$-unital $\A$. 
We prove that both definitions of the spectral flow coincide for paths where both apply.

We show that the spectral flow is uniquely determined by a set of natural axioms, which differ from axioms characterizing the spectral flow on a Hilbert space in that they include functoriality in $\A$. 

In analogy we develop the spectral flow  with values in $K_1(\A)$ for paths of odd selfadjoint regular Fredholm operators on the $\bbbz/2$-graded Hilbert $\A$-module $H_{\A}^+ \oplus H_{\A}^-$. Leichtnam--Piazza defined an odd spectral flow for paths of Dirac operators over $C^*$-algebras on even dimensional manifolds using $Cl(1)$-spectral sections (we call them Lagrangian spectral sections) \cite{lp2}. We give a new definition of the difference element of Lagrangian spectral sections which allows us to develop the theory in the odd case in analogy to the even case. Again, the odd relative index, which we define, is related to the generalized Maslov index in \cite{n}.

Though the odd spectral flow is trivial for a Hilbert space since $K_1(\bbbc)=0$, examples arise in classical analysis: Melrose's divisor flow \cite{m} is an odd spectral flow with $\A=C_0(\bbbr)$  via the identification $K_1(C_0(\bbbr))\cong \bbbz$. By relating the odd spectral flow to the spectral flow via suspension we obtain a $K$-theoretic proof of a formula expressing the spectral flow in terms of the divisor flow \cite{lmp}.

We show that under certain conditions the spectral flow can be expressed as the index of the operator $\ra_t +D_t$ and prove a corresponding result for the odd spectral flow. This is well-known for the classical spectral flow. For $\A \neq \bbbc$ results in this direction exist only for loops \cite{dz2}\cite{lp2}. The result makes it possible to use index theorems for computation of the spectral flow. 

One of the applications of these results is the definition of a noncommutative version of the Maslov index for a pair of paths of Lagrangians. The noncommutative Maslov index occurs as a correction term in a splitting formula of the spectral flow for loops of families of Dirac operators resp. Dirac operators over a $C^*$-algebra. A far-reaching formulation and proof of the classical splitting theorem for paths of Dirac type operators has been given by Nicolaescu \cite{n}. We derive the splitting formula of the spectral flow from the splitting theorem for the index. The latter has been proven in similar situations for example in \cite{n}\cite{dz1}\cite{lp2}. The proof sketched here might be of independent interest since, in contrast to the proofs in \cite{dz1}\cite{lp2}, it treats the odd case essentially on an equal footing with the even case. 

We relate the noncommutative Maslov index for paths to Bott periodicity and to the noncommutative  Maslov triple index  \cite{bk}\cite{wa}, thus generalizing properties of the classical Maslov index (see \cite{clm}). 

There is a generalization of the spectral flow to Breuer-Fredholm operators (see \cite{bcp} for an overview). The relation to our notion should be seen in the broader picture of the relation between index theory of Breuer-Fredholm operators and index theory of Fredholm operators over $C^*$-algebras: These are different frameworks with common applications, for example in higher index theory. 

Some results of this paper concerning the spectral flow on a separable Hilbert space are reconsidered in a more classical language in \cite{wa3}.

{\it Acknowledgements:} I would like to thank Peter Haskell for many helpful discussions, Paolo Piazza for providing me  Wu's unpublished preprint, and John Phillips and Max Karoubi for giving me some useful information about their work.

\tableofcontents

\section{Unbounded Fredholm operators on Hilbert $C^*$-modules}

Let $\A$ be a $\sigma$-unital $C^*$-algebra with norm $\| \cdot \|_{\A}$.

In this section we review some facts about Hilbert $C^*$-modules. More information on the subject can be found in \cite{la} and \cite{wo}. Important for the following are unbounded Fredholm operators on Hilbert $C^*$-modules. 

A Hilbert $\A$-module is a right $\A$-module $H$ with a map
  $\langle~,~\rangle:H \times H \to \A$, called an $\A$-valued scalar product, such that
\begin{enumerate}
\item $\langle~,~\rangle$ is $\A$-linear in the second variable,
\item $\langle x,y\rangle=\langle y,x\rangle^*$,
\item $\langle x,x\rangle \ge 0$ for all $x \in H$,
\item if $\langle x,x\rangle=0$, then $x=0$.
\end{enumerate}
Furthermore $H$ is assumed to be complete with respect to the norm $\|x\|:=\|\langle x,x\rangle \|_{\A}^{\frac 12}$.

If $\pi:M  \to B$ is a fiber bundle of closed Riemannian manifolds with compact base space, the closure of $C(M)$ with respect to norm induced by the $C(B)$-valued scalar product $\langle f,g\rangle(b)= \int_{M_b}\ov{f} g ~dvol_b$ is a Hilbert $C(B)$-module. By generalizing this example to vector bundles our results can be applied to family index theory. 

The standard Hilbert $\A$-module $H_{\A}$ is the right $\A$-module $\{(a_n)_{n \in \bbbn}~|~\sum\limits_{n
=1}^{\infty}a_n^*a_n \mbox{ converges}\}$ endowed with the $\A$-valued scalar
product $$\langle(a_n)_{n \in \bbbn},(b_n)_{n \in \bbbn}\rangle:=\sum_{n=1}^{\infty}
a_n^*b_n \ .$$ 

Most constructions will be done on $H_{\A}$. In order to generalize them to a countably generated Hilbert $\A$-module $H$ one can use the Stabilization Theorem, which states that the Hilbert $\A$-module $H \oplus H_{\A}$ is isomorphic to $H_{\A}$. One shows that the generalization does not depend on the choice of the isomorphism using that the unitary group ${\cal U}(H_{\A})$ is contractible. We often tacitly identify $H_{\A}^n$ with $H_{\A}$.

If $H_1,H_2$ are Hilbert $\A$-modules, then $B(H_1,H_2)$ is the space of bounded adjointable operators from $H_1$ to $H_2$. The subspace of compact operators, which is defined as the closure in $B(H_1,H_2)$ of the space spanned by the operators $\theta_{z,y}: x\mapsto z\langle y,x\rangle$ with $y\in H_1,z \in H_2,$ is denoted by $K(H_1,H_2)$. 

An operator $F\in B(H_{\A})$ is called Fredholm if its image $\pi(F)$ in $B(H_{\A})/K(H_{\A})$ is invertible. Its index $\ind(F) \in K_0(\A)$ is the image of $[\pi(F)] \in K_1(B(H_{\A})/K(H_{\A}))$ under the standard isomorphisms $K_1(B(H_{\A})/K(H_{\A})) \cong K_0(K(H_{\A})) \cong K_0(\A)$. 

A closed operator $D$ on a Hilbert $\A$-module $H$ is called regular if it is densely defined with densely defined adjoint and if $(1+D^*D)$ has a bounded inverse. There is a continuous functional calculus for selfadjoint regular operators, namely a unital $C^*$-algebra homomorphism $C(\bbbr) \to B(H),~f \mapsto f(D)$. The operator $D(1+D^*D)^{-\frac 12}$ is called the bounded transform of $D$. We say that a regular selfadjoint operator has compact resolvents if $(D + i)^{-1}$ and $(D - i)^{-1}$ are compact. 

For a regular operator $D$ we denote by $H(D)$ the Hilbert $\A$-module whose underlying vector space is $\dom D$ and whose scalar product is given by $$\langle v,w\rangle_{H(D)}=\langle v,w\rangle+ \langle Dv,Dw\rangle \ .$$
We say that a regular operator $D$ on $H_{\A}$ is Fredholm if it is Fredholm as a bounded operator from $H(D)$ to $H_{\A}$. This is equivalent to $D(1+D^*D)^{-\frac 12}$ being Fredholm on $H_{\A}$. 

We denote the space of selfadjoint regular Fredholm operators on $H_{\A}$ by $RSF(H_{\A})$ and the subspace of bounded selfadjoint Fredholm operators by $BSF(H_{\A})$.
  
The following elementary criteria for a selfadjoint regular operator to be Fredholm are  essential for the definition of the spectral flow: 

\begin{prop}
Let $D$ be a selfadjoint regular operator on $H_{\A}$. Then $D$ is Fredholm if and only if there is $\ve>0$ such that for all $\phi\in C(\bbbr)$ with $\supp \phi \subset [-\ve,\ve]$ we have that $\phi(D) \in K(H_{\A})$.
\end{prop}

\begin{proof}
The operator $D$ is Fredholm if and only if $D(1+D^2)^{-1/2}$ is Fredholm, hence if and only if the class of $D(1+D^2)^{-1/2}$ is invertible in $B(H_{\A})/K(H_{\A})$. Thus there is $\ve>0$ such that for all $\phi\in C(\bbbr)$ with support in $[-\ve,\ve]$ we have that $\phi(D(1+D^2)^{-1/2})=0$ in $B(H_{\A})/K(H_{\A})$. 
\end{proof} 

\begin{cor}
\label{critfred}
Let $D$ be a selfadjoint regular operator on $H_{\A}$. Then $D$ is Fredholm if and only if there is $\ve>0$ such that for any $\chi \in C(\bbbr)$ with $\chi|_{(-\infty,-\ve]}=-1$ and $\chi|_{[\ve,\infty)}=1$ we have that $\chi^2(D)-1 \in K(H_{\A})$.
\end{cor}

\begin{ddd}
If $\chi$ is an odd non-decreasing smooth function with $\chi(0)=0,~ \chi'(0)>0$, furthermore $\lim_{x \to \infty}\chi(x)= 1$ and $\chi(D)^2-1 \in K(H_{\A})$, then we call $\chi$ a {\rm normalizing function} for $D$.
\end{ddd}

For example if $D$ has compact resolvents, then $\chi(x)=x(x^2+1)^{-\frac 12}$ is a normalizing function for $D$. 

The notion ``normalizing function'' appears in \cite{hr} with a different but related meaning.

\subsection{Kasparov modules for unbounded Fredholm operators}

Index theory of Fredholm operators on Hilbert $C^*$-modules is organized in Kasparov's $KK$-theory. In order to assign a Kasparov $(\bbbc,\A)$-module to an element of $RSF(H_{\A})$ we introduce the notion of a truly unbounded Kasparov module. A reference for $KK$-theory is \cite{bl}. 

We fix the following conventions:
In a $\bbbz/2$-graded context $[~,~]$ denotes the supercommutator and $\ten$ the graded tensor product. If $V$ is a vector space, then $V^+$ resp. $V^-$ means the space $V$ endowed with the $\bbbz/2$-grading for which all elements are homogeneous of even resp. odd degree. If no grading is specified, we assume the grading for which all elements are even. We write $\hat H_{\A}= H_{\A}^+ \oplus H_{\A}^-$.

First recall the definition of an even resp. odd bounded Kasparov $(\B,\A)$-module. As before $\A$ is $\sigma$-unital.  The $C^*$-algebra $\B$ may be $\bbbz/2$-graded; here we have in mind the Clifford algebra $C_1$ with one odd generator. 
For simplicity we assume that $\B$ is unital.

A bounded {\it odd} resp. {\it even} Kasparov $(\B,\A)$-module is a triple $(H,\rho,T)$, where $H$ is a $\bbbz/2$-graded countably generated Hilbert $\A$-module, $\rho:\B \to B(H)$ is an even unital $C^*$-homomorphism and $T \in B(H)$ is a selfadjoint {\it even} resp. {\it odd} operator  with $T^2-1 \in K(H)$ and $[\rho(b),T] \in K(H)$ for all $b \in \B$. Its class in $KK_1(\B,\A)$ resp. $KK_0(\B,\A)$ is denoted by $[(H,\rho,T)]$ or simply $[T]$.

\begin{ddd}
Let $H$ be possibly $\bbbz/2$-graded countably generated Hilbert $\A$-module and let $\rho:\B \to B(H)$ be an even unital $C^*$-homomorphism. 

Let $D\in RSF(H)$ be even resp. odd and assume that there is a dense subset $\B' \subset \B$ such that for all $b \in \B'$ the operator $[D,\rho(b)]$ is defined on a core for $D$ and extends to a compact operator from $H(D)$ to $H$ and that there is $x \in [0,\frac 12)$ such that $[D,\rho(b)](1+D^2)^{-x}$ is bounded.

Then we call $(H,\rho,D)$ a truly unbounded odd resp. even Kasparov $(\B,\A)$-module. 

If $(H,\rho,D)$ is a truly unbounded odd resp. even $(\B,\A)$-Kasparov module and $\chi$ is a normalizing function for $D$, then $(H,\rho,\chi(D))$ is a bounded odd resp. even $(\B,\A)$-Kasparov module. We define the class $[(H,\rho,D)]$ in $KK_1(\B,\A)$ resp. $KK_0(\B,\A)$ to be the class $[(H,\rho,\chi(D))]$.
\end{ddd}  

\begin{proof}
We have to show that $[\chi(D),\rho(b)]$ is compact for all $b \in \B$. Since the algebra generated by $(x+i)^{-1}$ and $(x-i)^{-1}$ is dense in $C_0(\bbbr)$ and by continuity it is enough to prove that
$[D(D^2+1)^{-\frac 12},\rho(b)]$, $[(D+i)^{-1},\rho(b)]$ and $[(D-i)^{-1},\rho(b)]$ are compact for all $b \in \B'$. 
 
For $b \in \B'$ we have that
$$[D(D^2+1)^{-\frac 12},\rho(b)]=[D,\rho(b)](D^2+1)^{-\frac 12}+ D[(D^2+1)^{-\frac 12},\rho(b)] \ .$$
The operator $[D,\rho(b)](D^2+1)^{-\frac 12}$ is compact. For $f \in \dom D$
$$D[(D^2+1)^{-\frac 12},\rho(b)]f =\frac{1}{\pi} \int_0^{\infty}\lambda^{-\frac 12}D[(D^2+1+\lambda)^{-1},\rho(b)]f ~d\lambda \ .$$
Furthermore
\begin{eqnarray*}
D[(D^2+1+\lambda)^{-1},\rho(b)]f&=&D(D+i\sqrt{1+\lambda})^{-1}[(D-i\sqrt{1+\lambda})^{-1},\rho(b)]f\\
&& + D[(D+i\sqrt{1+\lambda})^{-1},\rho(b)](D-i\sqrt{1+\lambda})^{-1})f \ .
\end{eqnarray*}
Rewriting $$[(D\pm i\sqrt{1+\lambda})^{-1},\rho(b)]f=(D\pm i\sqrt{1+\lambda})^{-1}[D,\rho(b)](D\pm i\sqrt{1+\lambda})^{-1}f$$ we see that $\lambda^{-\frac 12}D[(D^2+1+\lambda)^{-1},\rho(b)]$ is compact and that the integral over $\lambda$ converges in the operator norm. It follows that  $[D(D^2+1)^{-\frac 12},\rho(b)]$ is compact. The compactness of $[(D \pm i)^{-1},\rho(b)]$ follows from the last equation.
\end{proof}

{\bf Remarks and examples:} 
(1) An elliptic selfadjoint pseudodifferential operator of non-negative order on a closed Riemannian manifold $M$ defines  a truly unbounded $(C(M),\bbbc)$-Kasparov module. 

(2) Truly unbounded Kasparov modules arise in connection with compactifications of complete manifolds: If the Dirac operator on a complete spin manifold $M$ is invertible near infinity, it defines a truly unbounded Kasparov $(\B,\bbbc)$-module, where $\B$ is the algebra of continuous functions on an appropriate compactification of $M$ (see \cite{bu} for details). 
\\

Since the group $KK_1(\bbbc,\A)$ is particularly important for the following we describe its construction more precisely:

We assume $\rho:\bbbc \to B(H_{\A})$ to be the unique unital $C^*$-homomorphism and suppress it in the notation.

An odd bounded Kasparov $(\bbbc,\A)$-module $(H_{\A},T)$ is called degenerate if $T^2=1$. Two bounded odd Kasparov $(\bbbc,\A)$-modules $(H_{\A},T_0)$ and $(H_{\A},T_1)$ are called homotopic if there is a bounded odd Kasparov $(\bbbc,C([0,1],\A))$-module $(H_{C([0,1],\A)},T)$ with $(\ev_0)_*T=T_0$ and $(\ev_1)_*T=T_1$, where $\ev_i:C([0,1],\A) \to \A$ denotes the evaluation at $i \in [0,1]$. As a set $KK_1(\bbbc,\A)$ is the quotient of the set of bounded odd Kasparov $(\bbbc,\A)$-modules $(H_{\A},T)$ by the equivalence relation generated by homotopy and direct sum with degenerate elements. With the direct sum as group operation $KK_1(\bbbc,\A)$ is an abelian group (here we use the identification of $H_{\A}^2$ with $H_{\A}$). 

The unitalization of a $C^*$-algebra $\A$ is denoted by $\A^{\sim}$. We have that $KK_1(\bbbc,\A)=KK_1(\bbbc,\A^{\sim})$.

The standard isomorphism $KK_0(\bbbc,\A) \cong K_0(\A)$ assigns to the class represented by an odd Fredholm operator $D=\left(\begin{array}{cc} 0 & D^- \\ D^+ & 0 \end{array}\right)$ on $\hat H_{\A}$ the index of $D^+$, whereas the standard isomorphism $KK_1(\bbbc,\A) \cong K_1(\A)$ assigns to an element in $KK_1(\bbbc,\A)$ represented by a bounded odd Kasparov module $(H_{\A},T)$ the element $$\ind^o(T):=[e^{\pi i (T+1)}] \in K_1(K(H_{\A})) \cong K_1(\A) \ .$$ If $D \in RSF(H_{\A})$ and $\chi$ is a normalizing function for $D$, then we write $$\ind^o(D):=\ind^o(\chi(D)) \in K_1(\A) \ .$$

\subsection{Families}
\label{fam}

Let $B$ be a compact space. 

We endow the space $C(B,H_{\A})$ with the obvious Hilbert $C(B,\A)$-module structure. Then naturally $C(B,H_{\A}) \cong H_{C(B,\A)}$. In this section we study under which conditions a family of regular selfadjoint operators on $H_{\A}$ with parameter space $B$ defines an element of $RSF(C(B,H_{\A}))$.

A family of bounded operators on $H_{\A}$ with parameter space $B$ defines an element in $B(C(B,H_{\A}))$ if and only if it and its adjoint are strongly continuous.

Let $(D_b)_{b \in B}$ be a family of regular selfadjoint operators on $H_{\A}$. It defines an operator on the $C(B,\A)$-module $C(B,H_{\A})$  with
$$\dom (D_b)_{b \in B}= \{x \in C(B,H_{\A})~|~x(b) \in \dom D_b;~ (b \mapsto D_bx(b)) \in C(B,H_{\A}) \} \ .$$
We write $(D_b)_{b \in B}$ or -- by abuse of notation -- simply $D_b$ for this operator on $C(B,H_{\A})$.  

For $b \in B$ we denote by $\ev_b:C(B,H_{\A}) \to H_{\A}$ the evaluation at $b$.

\begin{prop}
\label{critreg}
Let $B$ be a compact space.

Then $(D_b)_{b \in B}$ is a regular selfadjoint operator on $C(B,H_{\A})$ if and only if the resolvents $(D_b \pm i)^{-1}$ depend on $b$ in a strongly continuous way.
\end{prop}

\begin{proof}
If the right hand side holds, then it is clear that the family $(1+D_b^2)^{-1}$ defines a bounded operator on $C(B,H_{\A})$. Hence we only have to show that $\dom (D_b)_{b \in B}$, or equivalently $(D_b +i)^{-1}C(B,H_{\A})$ is dense in $C(B,H_{\A})$. 

For that end note first that $(D_b +i)^{-1}C(B,H_{\A})$ is a $C(B,\A)$-module. Furthermore for any $\ve >0$, $x \in H_{\A}$ and $b_0 \in B$ there is $f \in \dom (D_b)_{b \in B}$ with $\|f(b_0)-x\| < \ve$, namely let $y \in \dom D_{b_0}$ with $\|x-y\| < \ve$ and define $f(b)=(D_b +i)^{-1}(D_{b_0}+i)y$. By the following arguments these two properties imply that $\dom (D_b)_{b \in B}$ is dense in $C(B,H_{\A})$: 

Let $f \in C(B,H_{\A})$ and $\ve>0$. For each $b_0 \in B$ and $g \in \dom (D_b)_{b \in B}$ with $\|g(b_0)-f(b_0)\| < \ve/3$ let 
$$U_g(b_0):=\{b \in B~|~\|g(b)-g(b_0)\| < \ve/3,~\|f(b)-f(b_0)\| < \ve/3\} \ .$$
These sets define an open covering of $B$.
By compactness there is a finite covering $(U_{g_i}(b_i))_{i \in I}$ of $B$.
Let $(\chi_i)_{i \in I}$ be a partition of unity subordinate to the covering. Then $\sum_{i \in I} \chi_i g_i \in \dom (D_b)_{b \in B}$ and ´for all $b \in B$
 $$\|f(b)- \sum_{i \in I} \chi_i(b) g_i(b)\| \le \sum_{i \in I} \chi_i(b) \|f(b)-g_i(b)\| < \ve \ .$$ 
\end{proof}

In applications (cf. \S \ref{appl}) the following criterion is easier to verify.

\begin{lem}
\label{critreg2}
Let $B$ be a compact space.
Let $S:=\dom (D_b)_{b \in B}$.

If for each $b \in B$ the set $S_b:=\ev_b(S) \subset H_{\A}$ is a core for $D_b$, then $(D_b)_{b \in B}$ is a regular selfadjoint operator on $C(B,H_{\A})$. 
\end{lem}

\begin{proof}
The operator $(D_b+i)_{b \in B}$ is defined on $S$. Its closure on the Hilbert $C(B,\A)$-module $C(B,H_{\A})$ is injective and has closed range by \cite[Lemma 9.7]{la}. Since $S_b$ is a core for $(D_b+i)$ we have that $(D_b+i)S_b$ is dense in $H_{\A}$ for each $b \in B$. Furthermore $(D_b+i)_{b \in B}S$ is a $C(B,\A)$-module. By the argument in the previous proof this implies that $(D_b+i)_{b \in B}S$ is dense in $C(B,H_{\A})$, hence $(D_b +i)_{b \in B}$ is invertible on $C(B,H_{\A})$. Analogous arguments apply to $(D_b-i)_{b \in B}$. This implies that $(D_b)_{b \in B}$ is regular as an operator on $C(B,H_{\A})$.
\end{proof}

The following criterion helps to decide whether a family of Fredholm operators $(D_b)_{b \in B}$ is a Fredholm operator on $C(B,H_{\A})$.

\begin{lem}
\label{critrsf}
Let $B$ be a compact space. Let $\pi:B(C(B,H_{\A})) \to B(C(B,H_{\A}))/K(C(B,H_{\A}))$ be the projection.

A family $(D_b)_{b \in B}$ of selfadjoint regular Fredholm operators defines a regular selfadjoint Fredholm operator on $C(B,H_{\A})$ if and only if $b \mapsto (D_b \pm i)^{-1}$  
is strongly continuous and $\|\pi((1+D_b^2)_{b\in B}^{-1})\| <1$. 
\end{lem}

\begin{proof}

Regularity on $C(B,H_{\A})$ follows from Prop. \ref{critreg}.

Let $F_b:=D_b(1+D_b^2)^{-\frac 12}$.

If $$\|\pi((1+D_b^2)_{b\in B}^{-1})\|=\|(1-\pi((F_b)_{b\in B})^2)\| < c < 1 \ ,$$  then $\pi(F_b^2)$ is uniformly bounded below by $1-c$, hence there is a normalizing function $\chi$ for the operator $(F_b)_{b\in B}$ on $C(B,H_{\A})$. We conclude from Cor. \ref{critfred} that $(F_b)_{b\in B}$ is Fredholm on $C(B,H_{\A})$. 
\end{proof}

We note for later use:

\begin{lem}
\label{strcon}
Let $B$ be a compact space.
Let $(D_b)_{b \in B}$ be a regular selfadjoint operators with bounded inverse on $C(B,H_{\A})$.
 
Then $1_{\ge 0}(D_b)$ depends in a strongly continuous way on $b$.
\end{lem}

\begin{proof}
By $\|D_b x\|\ge \|D_b^{-1}\|^{-1} \|x\|$ we may find $f \in C(\bbbr)$ such that $f(D_b)=1_{\ge 0}(D_b)$ for all $b \in B$. By the functional calculus for regular operators $(f(D_b))_{b \in B} \in B(C(B,H_{\A}))$, hence $b \mapsto f(D_b)$ is strongly continuous.  
\end{proof}  

{\bf Example.} The following example illustrates the assumptions of Lemma \ref{critrsf}.
Let $(f_n)_{n \in \bbbn}$ be a sequence of non-decreasing functions in $C([0,1])$ with values in $[0,1]$. Assume that $f_n^{-1}(0)=\frac 12-\frac 1n$ and $f_n(t)=-1$ for $t \in [0, \frac 12-\frac 2n]$ and $f_n(t)=1$ for $t \in [\frac 12,1]$. Let $$F_t:l^2(\bbbn)\to l^2(\bbbn),~e_n \mapsto f_n(t) e_n  \ ,$$ where $\{e_n\}_{n \in \bbbn}$ is an orthonormal basis of $l^2(\bbbn)$. The family $(F_t)_{t \in [0,1]}$ is selfadjoint and strongly continuous with norm bounded by $1$, hence it is an element of $B(C([0,1],l^2(\bbbn)))$. Furthermore $F_t$ is Fredholm for any $t \in [0,1]$, but $(F_t)_{t \in [0,1]}$ fails to be Fredholm on $C([0,1],l^2(\bbbn))$. Note that $$\sup_{t \in [0,1]} \|\pi((1+F_t^2)^{-1})\|_{B(l^2(\bbbn))/K(l^2(\bbbn))} = \frac 12$$ but that $$\|\pi((1+F_t^2)^{-1})\|_{B(C([0,1],l^2(\bbbn)))/K(C([0,1],l^2(\bbbn)))}=1 \ .$$ In an intuitive sense the spectral flow of the path is $\infty$.

{\bf Remarks.} \\
(1) The properties tested in Lemma \ref{critfred} and Lemma \ref{critreg2} are local in the following sense: If $\{U_i\}_{i \in I}$ is a finite open covering of $B$ and $(D_b)_{b \in \ov{U_i}}$ is regular resp. Fredholm on $C(\ov{U_i},H_{\A})$ for any $i \in I$, then $(D_b)_{b \in B}$ is regular resp. Fredholm on $C(B,H_{\A})$. 

(2) Lemma \ref{critreg} and Lemma \ref{critrsf} can be used to find topologies on $RSF(H_{\A})$ guaranteeing that any continuous map $B \to RSF(H_{\A})$ with $B$ compact defines an element in $RSF(C(B,H_{\A}))$. The spectral flow will be defined and homotopy invariant for paths that are continuous with respect to such a topology.\\
An example is the gap topology on $RSF(H_{\A})$. It is the weakest topology such that the maps
$$RSF(H_{\A}) \to B(H_{\A}),~ D \mapsto (D+i)^{-1}$$ and $$RSF(H_{\A}) \to B(H_{\A}),~ D \mapsto (D-i)^{-1}$$ are continuous. We denote $RSF(H_{\A})$ endowed with the gap topology by $RSF(H_{\A})_{gap}$. It is shown in \cite{jo} that $RSF(H_{\A})_{gap}$ is a representing space for the functor $B \mapsto K_1(C(B,\A))$ from the category of compact topological spaces to the category of abelian groups. The natural isomorphism $[B,RSF(H_{\A})_{gap}] \cong K_1(C(B,\A))$ is defined by sending a class $[(D_b)_{b \in B}] \in [B,RSF(H_{\A})_{gap}]$ to the element $[(D_b)_{b \in B}]$ in $KK_1(\bbbc,C(B,\A))$. \\
A different topology is discussed in \cite{wa3}.

\section{Generalized spectral sections and the spectral flow}

From now on we assume that $\A$ is unital.

\subsection{Spectral sections}
\label{specsec}

We recall the definition of the noncommutative spectral flow for selfadjoint operators with compact resolvents in terms of spectral sections. We refer to \cite{lp2} for more about spectral sections and the spectral flow for Dirac operators over $C^*$-algebras on closed manifolds.  The proofs cited below from \cite{lp1}\cite{lp2}\cite{lp3} are formulated for Dirac operators but they remain valid for general selfadjoint regular operators with compact resolvents. (In fact, these references cite \cite{wu}, where the proofs are carried out as needed here.)  

Let $D$ be a selfadjoint regular operator with compact resolvents on $H_{\A}$.

A spectral cut  is a smooth function $\chi:\bbbr \to [0,1]$ such that there are $a,b \in \bbbr,~a<b$ with $\chi|_{(-\infty,a]}=0$ and $\chi|_{[b,\infty)}=1$.

A projection $P \in B(H_{\A})$ is called a spectral section of $D$ if there are spectral cuts $\chi_1,\chi_2$  with  
$$\Ran \chi_1(D) \subset \Ran  P \subset \Ran \chi_2(D)  \ .$$

We say that $D$ admits spectral sections if for any spectral cut $\chi$ there is a spectral section $P$ with $\Ran P \subset \Ran \chi(D)$ and a spectral section $Q$ with $\Ran \chi(D) \subset \Ran Q$. 

If $D$ admits spectral sections and $P$ is a spectral section of $D$, then there is a compact selfadjoint operator $A$ such that $D+A$ is invertible and $P=1_{\ge 0}(D+A)$ (see \cite[Prop. 2.10]{lp1}). 

If $D$ admits spectral sections, then given a pair of spectral sections $(P,Q)$ of $D$ there is a spectral section $R$ with $RP=RQ=R$. There is a well-defined difference element 
$$[P-Q]:=[P-R]-[Q-R] \in K_0(\A) \ .$$

Let $(D_t)_{t \in [0,1]}$ be a family of selfadjoint regular operators with compact resolvents that is regular as an operator on the Hilbert $C([0,1],\A)$-module $C([0,1],H_{\A})$. Assume furthermore that $(D_t)_{t \in [0,1]}$ admits spectral sections on $C([0,1],H_{\A})$ and let $(Q_t)_{t \in [0,1]}$ be a spectral section of $(D_t)_{t \in [0,1]}$. If $P_0$ resp. $P_1$ is a spectral section of $D_0$ resp. $D_1$, then the noncommutative spectral flow
$$\spfl(D_t,P_0,P_1):=[Q_0-P_0]-[Q_1-P_1]$$ does not depend on the choice of $(Q_t)_{t \in [0,1]}$.
If $D_0=D_1$, then $\spfl(D_t):=\spfl(D_t,P_0,P_0)$ does not depend on $P_0$.

The existence problem of spectral sections is nontrivial. A necessary condition is that $[D]=0$ in $KK_1(\bbbc,\A)$. The following lemma yields a sufficient condition. See also \cite[\S 2]{lp2}. 

\begin{lem} 
\label{exspsec} 
\begin{enumerate}
\item There is a spectral section for $D$ if and only if there is a projection $P \in B(H_{\A})$ with $\pi(P)=\pi(\chi(D)) \in B(H_{\A})/K(H_{\A})$ for a spectral cut $\chi$.
\item The operator $D$ admits spectral sections if there is a spectral section $P$ of $D$ such that $\Ran P$ and $\Ran (1-P)$ have complemented submodules isomorphic to $H_{\A}$.
\end{enumerate}
\end{lem}

\begin{proof}
(1) See \cite[Theorem 3 and Lemma 5]{lp2}. 

(2) See \cite[Prop. 3.5]{lp3}. The assumption implies that  $\Ran P \cong H_{\A}$  by the Stabilization Theorem. Hence there is a sequence of compact projections $(P_n)_{n \in \bbbn}$ with $P_n \le P$ converging strongly to the identity on $PH_{\A}$. An analogous statement holds true for $\Ran(1-P)H_{\A}$. This is needed in the proof.
\end{proof}

In order to guarantee that enough spectral sections exist we will use the following construction:

Let $d$ be a closed invertible operator with compact resolvents on the Hilbert space $H_{\bbbc}$ such that $1_{\ge 0}(d)$ and $1-1_{\ge 0}(d)$ are infinite dimensional projections. By tensoring with the identity we obtain a closed invertible operator with compact resolvents on $H_{\A}\cong H_{\bbbc} \ten \A$, also denoted by $d$. Assume that there is a spectral section for $D$. Then the operator $D \oplus d$ on $H_{\A}^2$ admits  spectral sections.

\subsection{Relative index of projections}

\label{relindex}

The relative index of a pair of projections $(P,Q)$ on $H_{\A}$ with $P-Q \in K(H_{\A})$, which we define in the following, generalizes the difference element of spectral sections. On a Hilbert space it is a special case of the relative index of a Fredholm pair of projections (see \cite{ass} for its properties). 

\begin{ddd}
A {\rm relative index} of projections is a map $\ind$ assigning to a pair $(P,Q)$ of projections on $H_{\A}$ with $P-Q \in K(H_{\A})$ an element in $K_0(\A)$ and fulfilling the following axioms: 
\begin{itemize}
\item[(I)] (Additivity.) If $P,Q,R$ are projections with $P-Q$ and $Q-R$ compact, then 
$$\ind(P,R)=\ind(P,Q)+\ind(Q,R) \ .$$
\item[(II)] (Functoriality.) If $\phi:\A \to \B$ is a unital $C^*$-homomorphism, then $$\phi_*\ind(P,Q)=\ind(\phi_*P,\phi_*Q) \ .$$ 
\item[(III)] (Stabilization.) If $P,Q$ are projections with $P-Q$ compact, and $R$ is a projection orthogonal to $P$ and $Q$, then $$\ind(P+R,Q+R)=\ind(P,Q) \ .$$ 
\item[(IV)] (Normalization.) If $P$ is a compact projection and $Q$ is a projection orthogonal to $P$, then $\ind(P+Q,Q)=[P] \in K_0(\A)$.
\end{itemize}
\end{ddd}

It will be shown in Prop. \ref{unique} that the relative index of projections is unique.

A relative index has the following properties:
\begin{lem}
\label{propgensp}
\begin{enumerate}
\item Let $(P,Q)$ be a pair of projections on $H_{\A}$ with $P-Q \in K(H_{\A})$. Then
$$\ind(P,Q)=-\ind(Q,P) \ .$$
\item If $D$ is a selfadjoint operator with compact resolvents and if $P,Q$ are spectral sections of $D$, then $\ind(P,Q)$ is uniquely defined. It agrees with $[P-Q]$ if $D$ admits spectral sections.
\item (Homotopy invariance.) If $(P_t)_{t \in [0,1]}$ and $(Q_t)_{t \in [0,1]}$ are strongly continuous paths of projections with $P_t-Q_t \in K(H_{\A})$ for all $t \in [0,1]$, then $\ind(P_0,Q_0)= \ind(P_1,Q_1)$.
\item If $U\in {\cal U}(H_{\A})$, then $\ind(UPU^*,UQU^*)=\ind(P,Q)$.
\end{enumerate}
\end{lem}

\begin{proof}
(1) follows from Additivity and since $\ind(P,P)=0$ by (IV).

(2) If $D$ admits spectral sections,  it is enough to consider $P \ge Q$ by Additivity. In this case the assertion follows from Normalization. If $D$ does not admit spectral sections, then let $d$ be as in the end of \S \ref{specsec}. Then $D \oplus d$ admits spectral sections and $P \oplus 1_{\ge 0}(d)$, $Q \oplus 1_{\ge 0}(d)$ are spectral sections of $D \oplus d$. By the Stabilization Axiom
$$\ind(P \oplus 1_{\ge 0}(d),Q \oplus 1_{\ge 0}(d))= \ind(P,Q) \ .$$ The left hand side is uniquely defined by the first part of the proof.

(3) The families $P_t$ and $Q_t$ define elements in $B(C([0,1],H_{\A}))$,
hence an element $\ind(P_t,Q_t) \in K_0(C([0,1],\A))$. By functoriality $(\ev_i)_*\ind(P_t,Q_t)=\ind(P_i,Q_i) \in K_0(\A),~i=0,1$. Now the result follows from  $(\ev_0)_*=(\ev_1)_*:K_0(C([0,1],\A) \to K_0(\A)$.

(4) follows from homotopy invariance and the contractibility of ${\cal U}(H_{\A})$.
\end{proof}

{\bf Examples:}

Let $P,Q \in B(H_{\A})$ be selfadjoint projections with $P-Q$ compact. 

(1) Let $<P,K(H_{\A})>$ be the -- in general not unital -- $C^*$-subalgebra of $B(H_{\A})$ generated by $P$ and $K(H_{\A})$.  The inclusion $K(H_{\A}) \to <P,K(H_{\A})>$ induces an exact sequence
$$0 \to K_0(K(H_{\A})) \to K_0(<P,K(H_{\A})>) \to K_0(\bbbc) \ ,$$ where the last map is zero if the projection $P$ is compact. We identify $K_0(K(H_{\A}))$ with $K_0(\A)$ and define
the $\ind_1(P,Q)$ as the inverse image in $K_0(\A)$ of $[P]-[Q] \in K_0(<P,K(H_{\A})>)$.
It is easily checked that $\ind_1$ is a relative index.

(2) If $P,Q$ are projections whose difference is compact, the operator $QP:P(H_{\A}) \to Q(H_{\A})$ is Fredholm with parametrix $PQ$. We define $\ind_2(P,Q)$ as the index of $QP:P(H_{\A}) \to Q(H_{\A})$. Additivity follows from the fact that
\begin{eqnarray*}
\ind_2(P,Q)+\ind_2(Q,R) &=& \ind(RQP:P(H_{\A}) \to R(H_{\A})) \\
&=& \ind(R(Q-R)P+RP:P(H_{\A}) \to R(H_{\A})) \\ 
&=& \ind(RP:P(H_{\A}) \to R(H_{\A}))  \\
&=& \ind_2(P,R) \ .
\end{eqnarray*}

(3) Let $U:\bbbr \to S^1$ be a continuous function with $U(x)=1$ for $x \in \bbbr \setminus (-1,1)$ and such that the continuous extension $U:\ov{\bbbr} \to S^1$ has winding number one. Here $\ov{\bbbr}$ is the one-point compactification of $\bbbr$. Then $$U( t (2Q-1)+(1-t)(2P-1)) \in \U(C_0(\bbbr,K(H_{\A}))^{\sim}) \ .$$ Hence it defines an element in $K_1(C_0(\bbbr,K(H_{\A})))$. We define  
$\ind_3(P,Q) \in K_0(\A)$ as its image under composition of the Bott periodicity map with the isomorphism $K_0(K(H_{\A})) \cong K_0(\A)$.

\subsection{Generalized spectral sections}

\begin{ddd}
Let $D$ be a regular selfadjoint operator on $H_{\A}$.
We call a selfadjoint operator $A \in B(H_{\A})$ such that $A:H(D) \to H_{\A}$ is compact and such that $D+A$ has a bounded inverse  a {\rm trivializing operator} for $D$.

If $A$ is a trivializing operator for $D$, we call $1_{\ge 0}(D+A)$ a {\rm generalized spectral section} of $D$ and $D+A$ a {\rm trivialization} of $D$.
\end{ddd}

{\bf Remarks.}
(1) The existence of a trivializing operator for $D$ implies that $D \in RSF(H_{\A})$ and $[D]=0$ in $KK_1(\bbbc,\A)$. 

(2) If $D$ has compact resolvents, then any spectral section of $D$ is a generalized spectral section. 

(3) There is a notion of trivializing operator in \cite{lp2}, which motivates the above definition. The notion of generalized spectral sections in \cite{dz2} is related but different from the notion used here.  

Assume from now on that $D \in RSF(H_{\A})$.

For each $R>0$ let $\phi_R:\bbbr \to [0,1]$ be a smooth function with $\phi_R(x)=1$ for $|x|\le R$ and $\phi_R(x)=0$ for $|x|>R+1$.

\begin{lem}
\label{boundspec}
Assume that $D$ has a bounded inverse and let $B$ be a symmetric operator on $H_{\A}$ with $\dom D \subset \dom B$ and such that $B:H(D) \to H_{\A}$ is compact. Then there is a generalized spectral section for $D+B$. It follows that $D+B \in RSF(H_{\A})$.
\end{lem}

\begin{proof}
Since $\phi_R(D)$ converges strongly to the identity for $R \to \infty$, the operator $\phi_R(D)B\phi_R(D):H(D) \to H_{\A}$ converges to $B \in K(H(D),H_{\A})$ in norm. Hence for $R$ big enough $D+B-\phi_R(D)B\phi_R(D):H(D) \to H_{\A}$ is invertible. 
\end{proof}

By the following lemma trivializing operators can always be chosen compact and ``smoothing''; if $D$ has compact resolvents, they can be approximated by trivializing operators yielding spectral sections.

\begin{lem} 
\label{spsec}
Let $A$ be a trivializing operator for $D$. Then for  $R$ big enough $K:=\phi_R(D)A\phi_R(D)$ is a trivializing operator for $D$.  
If $D$ has compact resolvents, then $P=1_{\ge 0}(D+K)$ is a spectral section of $D$.  
\end{lem}

\begin{proof}
The first assertion holds since $\phi_R(D) A\phi_R(D)$  converges in norm to $A$ as an operator from $H(D)$ to $H_{\A}$ for $R \to \infty$.

Now assume that $D$ has compact resolvents and let $K:=\phi_R(D)A\phi_R(D)$ be a trivializing operator for $D$.

For $f \in C(\bbbr)$ with $\supp f \cap [-R-1,R+1]= \emptyset$ and $x \in \dom D$ we have that
$$f(D)(D+ K)x=(D+ K)f(D)x = Df(D)x\ ,$$
hence  $$f(D)(D+ K\pm i)^{-1}= (D+ K \pm i)^{-1}f(D)= (D\pm i)^{-1}f(D)\ .$$
It follows that for any $g \in C(\bbbr)$ for which $\lim\limits_{x\to \infty}g(x)$ and $\lim\limits_{x\to -\infty}g(x)$ exist
$$f(D)g(D+ K)=g(D+ K)f(D)=g(D)f(D) \ .$$
Hence if $f$ is a spectral cut with $\supp f \subset (R,\infty)$, then $$1_{\ge 0}(D+ K) f(D) =f(D) \ ,$$
thus $\Ran f(D) \subset \Ran 1_{\ge 0}(D+ K)$. If $f$ is a spectral cut with $\supp(1-f)\subset (-\infty,-R)$, then $$f(D)1_{\ge 0}(D+ K)=1_{\ge 0}(D+ K) \ ,$$ hence $$\Ran 1_{\ge 0}(D+ K) \subset \Ran f(D) \ .$$ This shows the assertion.
\end{proof}

Let $\ind$ be a relative index of projections.

\begin{prop}
\label{diffspec}
Let $D_0,D_1$ be regular selfadjoint operators such that $D_0-D_1:H(D_0) \to H_{\A}$ is compact and $D_0-D_1 \in B(H_{\A})$.  Then for any function $f \in C(\bbbr)$ such that $\lim_{x \to \infty} f(x)$ and $\lim_{x \to -\infty} f(x)$ exist, we have that $f(D_0)-f(D_1) \in K(H_{\A})$. 
\end{prop}

\begin{proof}
First assume that $D_0,D_1$ are bounded. Then $D_0-D_1 \in K(H_{\A})$, hence $\pi(D_0)=\pi(D_1)$, where $\pi:B(H_{\A}) \to B(H_{\A})/K(H_{\A})$ is the projection. Thus $\pi(f(D_0))=\pi(f(D_1))$. 

In the unbounded case we claim that
$$D_0(1+D_0^2)^{-1/2}-D_1(1+D_1^2)^{-1/2} \in K(H_{\A}) \ .$$
Then the assertion follows from the bounded case.

For $\lambda \ge 0$ write $R_k(\lambda)=(D_k+ i\sqrt{1+\lambda})^{-1}$.
Then $R_1(\lambda)-R_0(\lambda)=R_1(\lambda)(D_0-D_1)R_0(\lambda)$ and $R_1(\lambda)^*-R_0(\lambda)^*=R_1(\lambda)^*(D_0-D_1)R_0(\lambda)^*$, hence 
\begin{eqnarray*}
\lefteqn{(1+D_1^2+ \lambda)^{-1}-(1+D_0^2+ \lambda)^{-1}}\\
&=& R_1(\lambda)(R_1(\lambda)^*-R_0(\lambda)^*)+(R_1(\lambda)-R_0(\lambda))R_0(\lambda)^*\\
&=&(1+D_1^2+ \lambda)^{-1}(D_0-D_1)(D_0 - i\sqrt{1+\lambda})^{-1}\\
 &&+ (D_1 + i\sqrt{1+\lambda})^{-1}(D_0-D_1)(1+D_0^2+ \lambda)^{-1} \ .
\end{eqnarray*}
This shows that $(1+D_1^2+ \lambda)^{-1}-(1+D_0^2+ \lambda)^{-1} \in K(H_{\A},H(D))$ with norm bounded by $C(1+\lambda)^{-1}$ for some $C>0$.  It follows easily that $D_1(1+D_1^2+ \lambda)^{-1}-D_0(1+D_0^2+ \lambda)^{-1} \in K(H_{\A})$ with norm bounded by $C(1+\lambda)^{-1}$. The claim holds since for $f \in \dom D_0= \dom D_1$ 
$$D_k(1+D_k^2)^{-1/2}f= \frac{1}{\pi} \int_0^{\infty} \lambda^{-1/2}D_k(1+D_k^2+ \lambda)^{-1}f ~d \lambda \ .$$
\end{proof}

Note that the proof needs only a slight modification to work also if $(D_0-D_1)(D_0^2+1)^{-\frac 12} \in K(H_{\A})$ and $(D_0-D_1)(D_0^2+1)^{-x} \in B(H_{\A})$ for some $0\le x<\frac 12$. 

In particular if $A_0,A_1$ are trivializing operators for $D \in RSF(H_{\A})$ and $P_i:= 1_{\ge 0}(D+A_i)$, then $P_0-P_1 \in K(H_{\A})$.

We write $\ind(D,A_0,A_1):=\ind(P_1,P_0)$.

\begin{lem}
\label{uniqcom}
Let $D$ be a regular selfadjoint operator with compact resolvents.
Let $A_0,A_1$ be trivializing operators of $D$. 
Then $\ind(D,A_0,A_1)$ does not depend on the choice of the relative index of projections.
\end{lem}

\begin{proof}
By Lemma \ref{spsec} and the homotopy invariance of the relative index we may assume that $1_{\ge 0}(D+A_i)$ are spectral sections of $D$. Then the assertion follows from Lemma \ref{propgensp}.
\end{proof}

\begin{ddd}
Let $B$ be a compact space.
Let $(D_b)_{b \in B}$ be a regular operator on the Hilbert $C(B,\A)$-module $C(B,H_{\A})$. We call a family $(A_b)_{b \in B}$ of bounded operators that is a trivializing operator for $(D_b)_{b \in B}$ on $C(B,H_{\A})$ a {\rm trivializing family} for $(D_b)_{b \in B}$.

We say that locally there are trivializing families for $(D_b)_{b \in B}$ if for any point $b_0 \in B$ there is a compact neighborhood $U_{b_0} \ni b_0$ such that there exists a trivializing family for $(D_b)_{b \in U_{b_0}}$. 
\end{ddd}

\subsection{The spectral flow}
\label{spflone}

Now we give the definition of the spectral flow in terms of generalized spectral sections. Fix a relative index of projections $\ind$. We will denote it by $\spfl_i$ in order to distinguish it from the definition in \ref{spflBott}. 

\begin{ddd}
Let $(D_t)_{t \in [0,1]}$ be a selfadjoint regular operator on the Hilbert $C([0,1],\A)$-module $C([0,1],H_{\A})$ for which  locally trivializing families exist. Let
 $0=t_0<t_1< \dots <t_n=1$ be such that there is a trivializing family $(B^i_t)_{t \in [t_i,t_{i+1}]}$ of $(D_t)_{t \in [t_i,t_{i+1}]},~i=0,\dots, n-1$.

For trivializing operators $A_0$ of $D_0$ and $A_1$ of $D_1$ we define
\begin{eqnarray*}
\spfl_i((D_t)_{t \in [0,1]},A_0,A_1) &:=& \ind(D_0,A_0,B^0_0) + \ind(D_1,B^{n-1}_1,A_1)\\
&&+ \sum_{i=1}^{n-1} \ind(D_{t_i},B^{i-1}_{t_i},B^i_{t_{i}}) \in K_0(\A) \ .
\end{eqnarray*}
\end{ddd}

Using additivity, homotopy invariance of the relative index of projections and Lemma \ref{strcon} one can check that if there are two families $(B_t)_{t \in [a,b]},(C_t)_{t \in [a,b]}$ of trivializing operators for $(D_t)_{t \in [a,b]}$ for some $0 \le a <b \le 1$, then $\ind(D_a,B_a,C_a)=\ind(D_b,B_b,C_b)$. 
From this one deduces that the definition of the spectral flow is independent of the subdivision and the choice of the  trivializing families $(B^i_t)_{t \in [t_i,t_{i+1}]}$. 

The spectral flow has the following properties: 
\begin{enumerate}
\item {\it (Homotopy invariance.)} If $(D_{s,t})_{s,t \in [0,1]}$ is a regular operator on $C([0,1] \times [0,1],H_{\A})$ for which locally trivializing families exist and $(A_s)_{s\in [0,1]}$ resp. $(B_s)_{s \in [0,1]}$ is a trivializing family for $(D_{s,0})_{s \in [0,1]}$ resp. $(D_{s,1})_{s \in [0,1]}$, then $$\spfl_i((D_{0,t})_{t \in [0,1]},A_0,B_0)=\spfl_i((D_{1,t})_{t \in [0,1]},A_1,B_1) \ .$$ 
\item {\it (Concatenation.)}
If $(D_t)_{t \in [0,1]}, (E_t)_{t \in [0,1]}$ are selfadjoint regular operators on $C([0,1],H_{\A})$  with $D_1=E_0$, then we define the concatenation $((D*E)_t)_{t \in [0,1]} \in RSF(C([0,1],H_{\A}))$ by $(D*E)_t=D_{3t}$ for $t \in [0,\frac 13]$, $(D*E)_t=D_1=E_0$ for $t \in [\frac 13,\frac 23]$ and $(D*E)_t=E_{3t-2}$ for $t \in [\frac 23,1]$. If locally trivializing families exist for $(D_t)_{t \in [0,1]}$ and $(E_t)_{t \in [0,1]}$, then they exist also for the concatenation. For any trivializing operator $B$ of $D_1$ we have that
$$\spfl_i((D * E)_t,A_0,A_1)=\spfl_i(D_t,A_0,B) + \spfl_i(E_t,B,A_1) \ .$$
\item {\it (Loops.)} If $(D_t)_{t \in [0,1]}$ is a family with $D_0=D_1$ for which locally trivializing families exist, then the spectral flow $$\spfl_i((D_t)_{t \in [0,1]}):=\spfl_i((D_t)_{t \in [0,1]},A,A)$$ does not depend on the choice of the trivializing operator $A$.
\item If $(D_t)_{t \in [0,1]}$ is a family with $D_0,D_1$ invertible and for which locally trivializing families exist  and if $\chi$ is a normalizing function for $(D_t)_{t \in [0,1]} \in RSF(C([0,1],H_{\A}))$, then 
$$\spfl_i(D_t,0,0)=\spfl_i(\chi(D_t),0,0) \ .$$
\item If $P,Q$ are projections with $P-Q \in K(H_{\A})$, then
$$\ind(Q,P) = \spfl_i((1-t)(2P-1)+t(2Q-1),0,0) \ .$$
\end{enumerate}

\begin{prop}
\label{unique} 

The relative index of projections is uniquely defined. 
\end{prop}

\begin{proof} 
Note that by Lemma \ref{uniqcom} the spectral flow of a path of operators with compact resolvents, if defined, does not depend of the choice of the relative index.

Let $P,Q$ be projections on $H_{\A}$ with $P-Q \in K(H_{\A})$.

Let $K$ be a strictly positive compact operator with norm smaller than $1/4$. 
In \cite[proof of Prop. 2.13]{jo}, it was shown that if $F$ is a selfadjoint Fredholm operator with $\|F\| \le 1$ and $F^2-1 \in K(H)$, then $(1-K)F(1-K)$ is the bounded transform of a regular selfadjoint operator with compact resolvents.

Hence $(1-K)((1-t)(2P-1)+t(2Q-1))(1-K)$ is the bounded transform of a regular operator $(D_t)_{t\in [0,1]}$ on $C([0,1],H_{\A})$ with compact resolvents.

For any spectral cut $\chi$ the class of $(\chi(D_t))_{t \in [0,1]}$ in $B(C([0,1],H_{\A}))/K(C([0,1],H_{\A}))$ equals the class of the constant family $P$, hence by Lemma \ref{exspsec} there is a (global) spectral section for $(D_t)_{t \in [0,1]}$. Thus $\spfl_i((D_t)_{t\in [0,1]},0,0)$ is well-defined. By Lemma \ref{uniqcom} it is independent of the choice of the relative index. The assertion follows from 
\begin{eqnarray*}
\ind(Q,P) &=& \spfl_i((1-t)(2P-1)+t(2Q-1),0,0) \\
&=&\spfl_i((1-K)((1-t)(2P-1)+t(2Q-1))(1-K),0,0)\\
&=& \spfl_i((D_t)_{t\in [0,1]},0,0)\ .
\end{eqnarray*}
The first equality follows from Example (1) in \S \ref{spflone} and the last one from Property (4) of $\spfl_i$. 
\end{proof}

{\bf Remark.} 
 Let $(D_t)_{t \in [0,1]}$ be a gap continuous family of selfadjoint Fredholm operators on a separable Hilbert space. The spectral flow $\spfl_i$ is well-defined, since, by \cite[Prop. 2.10]{blp}, for each $t_0 \in [0,1]$ there is $a >0$ such that $1_{[-a,a]}(D_t)$ is compact and continuous in $t$ in a neighborhood of $t_0$. Thus $D_t + 2a 1_{[-a,a]}(D_t)$ is invertible near $t_0$, hence $2a 1_{[-a,a]}(D_t)$ is a trivializing family near $t_0$. 
It is easy to verify that the spectral flow defined in \cite{blp} for gap continuous paths agrees with our notion.

The following lemma shows that in favorable cases we can conclude from the pointwise existence of trivializing operators to the local existence of trivializing families.

\begin{lem}
Let $B$ be a compact space.
Let $(D_b)_{b \in B}$ be a family of regular selfadjoint operators with common domain. Assume that for each $b_0 \in B$ and $b \in B$ the operator $D_b:H(D_{b_0}) \to H_{\A}$ is continuous and depends continuously on $b$ with respect to the operator norm. 

If for some $b_0 \in B$ there is a trivializing operator $A$ for $D_{b_0}$, then there is a neighborhood $U \ni b_0$ such that $D_b+A$ is invertible for all $b \in U$. The inverse depends continuously on $b$. 
\end{lem}

\begin{proof}
For $$\|(D_b-D_{b_0})(D_{b_0}+A)^{-1})\| < \frac 12 $$ the Neumann series  
$$(D_b+A)^{-1}=\sum_{n=0}^{\infty} (D_{b_0}+A)^{-1}((D_b-D_{b_0})(D_{b_0}+A)^{-1})^n$$ converges and $$\|(D_b+A)^{-1}-(D_{b_0}+A)^{-1}\| \le 2\|(D_{b_0}+A)^{-1}\| \|(D_b-D_{b_0})(D_{b_0}+A)^{-1})\| \ .$$
\end{proof}

{\bf Example.}
Let $D \in RSF(H_{\A})$ be invertible and let $(B_t)_{t \in [0,1]}$ be a family of symmetric operators such that $\dom D \subset \dom B_t$ for all $t \in [0,1]$ and such that $B_t:H(D) \to H_{\A}$ is compact and depends continuously on $t$. Assume that $B_0=0$ and that $D+B_1$ is invertible. From Lemma \ref{boundspec} we conclude that $(D+B_t)_{t \in [0,1]} \in RSF(C([0,1],H_{\A}))$ and that $\spfl_i((D+B_t)_{t \in [0,1]},0,0)$ is well-defined. We prove that the spectral flow of $(D +B_t)_{t \in [0,1]}$ depends only on the endpoints. Let $\phi_R$ as defined before Lemma \ref{boundspec}. For $R$ big enough the operator $D + sB_1 + (1-s)\phi_R(D)B_1\phi_R(D)$ is invertible for all $s \in [0,1]$. Hence by homotopy invariance $\spfl_i(D + sB_t + (1-s)\phi_R(D)B_t\phi_R(D))$ is independent of $s$. Thus
\begin{eqnarray*}
\lefteqn{\spfl_i(D+B_t,0, 0)}\\
&=&\spfl_i(D +\phi_R(D)B_t\phi_R(D),0,0)\\
&=&\ind(1_{\ge 0}(D +\phi_R(D)B_1\phi_R(D)),1_{\ge 0}(D)) \ .
\end{eqnarray*} 
This shows the claim.
If $B_1(1+D^2)^{-x}$ is compact for some $0\le x<\frac 12$, then by the remark following Prop. \ref{diffspec} 
$$\spfl_i(D+B_t,0, 0)=\ind(1_{\ge 0}(D+B_1),1_{\ge 0}(D)) \ .$$
For separable Hilbert spaces this formula is shown in \cite[Theorem 3.6]{le} under the more general assumption that $D+B_t$ is Riesz continuous and without further restriction on $B_1$.

\subsection{Spectral flow and index}

Besides establishing some technical lemmata needed in the sequel we want to clarify the relation of the spectral flow with the index of $\ra_t +D_t$. We have to impose additional conditions on $D_t$, in particular that it has compact resolvents, in order to guarantee that $\ra_t +D_t$ is Fredholm.

In the following we extend a family $(D_t)_{t \in [0,1]}$ without further notice to $\bbbr$ by setting $D_t=D_0$ for $t \in (-\infty,0]$ and $D_t=D_1$ for $t \in [1,\infty)$.

\begin{lem}
\label{relind}
Let $(D_t)_{t \in [0,1]}$ be a family of selfadjoint operators with compact resolvents such that $$\bbbr \to B(H_{\A}),~ t \mapsto (D_t-D_0)$$ is well-defined and strongly continuous.  Then in particular $(D_t)_{t \in [0,1]}$ is in $RSF(C([0,1],H_{\A}))$. Assume furthermore that $D_0$ and $D_1$ are invertible. It holds:
\begin{enumerate}
\item The closure ${\cal D}$ of the operator
$\ra_t +  D_t$ with domain $\C_c(\bbbr, \dom D_0)$ is regular on the Hilbert $\A$-module $L^2(\bbbr) \ten H_{\A}$.
\item If $\psi \in {\cal S}(\bbbr)$, then multiplication by $\psi$ is a compact operator from $H({\cal D})$ to $L^2(\bbbr) \ten H_{\A}$.
\item ${\cal D}$ is Fredholm. 
\item If $D_t$ is invertible for any $t \in \bbbr$, then $\ind({\cal D})=0$. 
\item If $(D_t)_{t \in [0,1]}$ is the concatenation of two families $(E_t)_{t \in [0,1]}$ and $(F_t)_{t \in [0,1]}$ with $E_1=F_0$ invertible, then
$$\ind {\cal D}=\ind {\cal E} + \ind {\cal F} \in K_0(\A) \ ,$$ where ${\cal E},{\cal F}$ are defined as in (1).
\end{enumerate}
\end{lem}

In the following taking closures is tacitly assumed if necessary.

\begin{proof}
(1) The Fourier transform yields a unitary transformation on $L^2(\bbbr) \ten H_{\A}$ which intertwines $1-\ra_t^2+D_0^2$ and $(1+t^2+D_0^2)$. By the regularity of $D_0$ the inverse of $(1+t^2+D_0^2)$ exists as a bounded operator on $L^2(\bbbr) \ten H_{\A}$. Hence $\ra_t +D_0$ is regular.  
The regularity of ${\cal D}$ follows since ${\cal D}$ is a bounded perturbation of $\ra_t + D_0$. 

(2) As in (1) one shows that $(\ra_t+D_0)^{-1}$ is a bounded operator on $L^2(\bbbr) \ten H_{\A}$. The operator $\psi(\ra_t+D_0)^{-1}$ is conjugate to $$L^2(\bbbr) \ten H_{\A} \ni  f \mapsto \int_{\bbbr} \hat \psi(\cdot -y)(D_0-iy)^{-1}f(y) ~dy $$ via Fourier transform.  Hence $\psi(\ra_t+D_0)^{-1}$ is an integral operator with integral kernel in $L^2(\bbbr \times \bbbr,K(H_{\A})),$ thus it is in $K(L^2(\bbbr) \ten H_{\A})$.
It follows that $\psi:H({\cal D}) \to L^2(\bbbr) \ten H_{\A}$ is a compact operator.

(3) Let $\phi \in \C(\bbbr)$ with $\phi((-\infty,0])=1$ and $\phi([1,\infty))=0$ and let $\chi_1 \in \C(\bbbr)$ with $\chi_1((-\infty,2])=1$ and $\chi_1([3,\infty))=0$ and $\chi_2 \in \C(\bbbr)$ with $\chi_2((-\infty,-2])=0$ and $\chi_2([-1,\infty))=1$.
Using (2) one easily verifies that $$R:= \phi(\ra_t +D_0)^{-1}\chi_1  + (1-\phi) (\ra_t +D_1)^{-1}\chi_2$$ 
is a parametrix for ${\cal D}$. 

(4) For $s \in [0,1]$ let $D(s)_t=D_t$ for $t \in (-\infty,s]$ and $D(s)_t=D_s$ for $t \in [s,\infty)$.  By part (1) the closure of $\ra_t+ D(s)_t$ is a continuous path of Fredholm operators from $H({\cal D})$ to $L^2(\bbbr) \ten H$ connecting ${\cal D}$ with the closure of $\ra_t + D_1$. Since the latter is invertible, we have that  $\ind({\cal D})=0.$

(5) The proof is analogous to the proof of the relative $K$-theoretic index theorem in \cite{bu}. We give it here for completeness.

Let ${\cal H}:=L^2(\bbbr,H_{\A}^+) \oplus L^2(\bbbr,H_{\A}^-)$. Let ${\cal G}_1^+= {\cal D}$ as an operator from $L^2(\bbbr,H_{\A}^+)$ to $L^2(\bbbr,H_{\A}^-)$, let ${\cal G}_1^- = {\cal D}^*$ and let
$${\cal G}_1=\left(\begin{array}{cc} 0 & {\cal G}_1^-\\ {\cal G}_1^+ & 0 \end{array}\right)$$ as an operator on ${\cal H}$. Analogously we define ${\cal G}_2,{\cal G}_3,{\cal G}_4$ with ${\cal G}_2^+=\ra_t+E_1=\ra_t+F_0$, ${\cal G}_3^+=\ra_t +(E*E_1)_t$ and ${\cal G}_4^+=\ra_t+ (F_0* F)_t$. Here $E_1$ and $F_0$ are considered as constant families. Define $\tilde {\cal H}:={\cal H}^2 \oplus ({\cal H}^{op})^2$ where ${\cal H}^{op}$ denotes the Hilbert $\A$-module ${\cal H}$ with reversed grading. Then ${\cal G}={\cal G}_1\oplus {\cal G}_2 \oplus {\cal G}_3 \oplus {\cal G}_4$ is an odd selfadjoint operator on $\tilde {\cal H}$ defining an element $[\tilde {\cal H},{\cal G}]$ in $KK_0(\bbbc,\A)$. We have that in $KK_0(\bbbc,\A)\cong K_0(\A)$
$$[\tilde {\cal H},{\cal G}]=\ind {\cal D}-\ind {\cal E} - \ind {\cal F} \ .$$  
Let $\chi_1,\chi_2 \in \C(\bbbr)$ with $\chi_1^2+ \chi_2^2=1$  and with $\chi_1((-\infty, \frac 13])=\chi_2([\frac 23, \infty))=1$. Let $z$ be the grading operator. Then the operator
$$X=z\left(\begin{array}{cccc} 0& 0& -\chi_1& -\chi_2 \\
0&0 & -\chi_2 & \chi_1 \\
\chi_1 & \chi_2 & 0 & 0  \\
\chi_2 & -\chi_1 & 0 & 0 
\end{array}\right) \ .$$ 
is an odd selfadjoint involution on $\tilde {\cal H}$. One verifies easily that $X{\cal G}+{\cal G}X$ is bounded and has compact $\bbbr$-support.
Hence
$[\tilde {\cal H},{\cal G}] \in \Ima(i^*:KK_0(C_1,\A) \to KK_0(\bbbc,\A))$. Since $i^*$ equals the Kasparov product from the left with $[i] \in  KK_0(\bbbc,C_1) =0$, it follows that $[\tilde{\cal H},{\cal G}]=0$ in $KK_0(\bbbc,\A)$.
\end{proof}

In concrete situations the conditions on the path can be relaxed. The crucial point is the Fredholm property which is guaranteed for example for smooth paths of elliptic pseudodifferential operators of order one on a closed manifold.

\begin{prop}
\label{indspecsec}
Let $D$ be a regular selfadjoint operator with compact resolvents
and let $A_0,A_1$ be trivializing operators of $D$. 
Let $\chi_i \in \C(\bbbr),~i=0,1,$ with $\chi_0|_{(-\infty,-1]}=1$ and $\chi_0|_{[0,\infty)}=0$, furthermore $\chi_1|_{[2,\infty)}=1$ and $\chi_1|_{(-\infty,1]}=0$.
Then 
$$\ind(\ra_t + D+\chi_0(t) A_0 + \chi_1(t) A_1)=\ind(D,A_0,A_1) \ .$$ 
\end{prop}

\begin{proof}
The operator is Fredholm by Lemma \ref{relind}. By Lemma \ref{spsec} we can reduce to the case where $A_0$ and $A_1$ are associated to spectral sections $P_0,P_1$. 

We first assume that $D$ admits spectral sections and treat the general case below. 
By Lemma \ref{relind} (5) and the properties of spectral sections we may reduce to the case $P_0\le P_1$. 

Let $R,Q$ be spectral sections with $Q\le  P_0 \le  P_1 \le R$ and such that $QDQ$ resp. $(1-R)D(1-R)$ are invertible on $Q(H_{\A})$ resp. $(1-R)(H_{\A})$.

If $A$ and $B$ are trivializing operators with $1_{\ge 0}(D+A)=1_{\ge 0}(D+B)$, then $D+A$ and $D+B$ are homotopic via the path $D_t:=(1-t)(D+A) + t(D+B)$. Furthermore $D_t=((1-t)|D+A| + t|D+B|)(2P-1)$ is invertible for any $t \in [0,1]$. 
Hence by Lemma \ref{relind} (4) and (5) we may assume that $$D+A_i=QDQ \oplus (P_i-Q) \oplus (P_i-R) \oplus (1-R)D(1-R) \ .$$
Then
\begin{eqnarray*}
\lefteqn{\ind(\ra_t +D + \chi_0(t) A_0+ \chi_1(t) A_1)}\\
&=& \ind\bigl((\ra_t + \chi_0(t) +\chi_1(t))|_{L^2(\bbbr,(P_0-Q)H_{\A})}\bigr)\\
&& +\ind\bigl((\ra_t -\chi_0(t) + \chi_1(t))|_{L^2(\bbbr,(P_1-P_0)H_{\A})}\bigr)\\
&&  + \ind\bigl((\ra_t -\chi_0(t) -\chi_1(t))|_{L^2(\bbbr,(R-P_1)H_{\A})}\bigr)\\
&=&\ind\bigl((\ra_t -\chi_0(t) + \chi_1(t))|_{L^2(\bbbr,(P_1-P_0)H_{\A})}\bigr) \\
&=&[P_1-P_0] \ .
\end{eqnarray*} 

Now the general case:

Let $d$ be as in the end of \S \ref{specsec}. Let $\tilde P_i$ be the spectral sections $1_{\ge 0}(d\oplus (D+A_i))$ on $H_{\A}^2$. Since $d \oplus D$ admits spectral sections, we get from the first part of the proof
$$\ind(\ra_t +d\oplus (D + \chi_0(t) A_0+ \chi_1(t) A_1))=\ind(\tilde P_1,\tilde P_0) \ .$$
This equation is equivalent to the assertion since Lemma \ref{relind} (4) implies that
$\ind(\ra_t +d)=0$ and since furthermore $\ind(\tilde P_1,\tilde P_0)=\ind(P_1,P_0)$.  
\end{proof}

\begin{prop} 
\label{spflind}
Let $(D_t)_{t \in [0,1]}$ be a family of regular selfadjoint operators with compact resolvents such that $D_t-D_0$ is bounded and depends in a strongly continuous way of $t$.

Let $A_0,A_1$ be the trivializing operators for $D_0,D_1$ and let $\chi_0,\chi_1$ be as in the previous proposition.
Then  
$$\spfl_i((D_t)_{t \in [0,1]},A_0,A_1)=\ind(\ra_t +D_t + \chi_0(t) A_0 + \chi_1(t) A_1) \in K_0(\A) \ .$$
\end{prop}

\begin{proof}
Let ${\cal D}=\ra_t +D_t + \chi_0(t) A_0 + \chi_1(t) A_1$.

Let  $B_t$ be a trivializing family for $D_t$ (which exists by Lemma \ref{boundspec}). 
Let $\phi \in \C(\bbbr)$ be such that $\phi|_{[0,1]}=1$ and $\supp \phi \in [-\frac 13,\frac 43]$. 
By Lemma \ref{relind} (2) multiplication by $\phi$ is a compact operator from $H({\cal D})$ to $L^2(\bbbr) \ten H_{\A}$, hence the index of ${\cal D}$ equals the index of 
$$\ra_t + D_t + \chi_0(t)A_0 +  \phi(t)B_t  + \chi_1(t)A_1 \ .$$ 
Thus by Lemma \ref{relind} (5)
\begin{eqnarray*}
\ind {\cal D} &=& \ind\bigl(\ra_t + D_0 + \chi_0(t)A_0 + \chi_1(t)B_0\bigr) \\
&&+ \ind\bigl(\ra_t + D_t  +B_t\bigr)\\
&& + \ind\bigl(\ra_t + D_1 + \chi_0(t) B_1  + \chi_1(t)A_1\bigr) \ .
\end{eqnarray*} 
Now the assertion follows from the previous lemmata.
\end{proof}

{\bf Remarks.}
(1) In special situations the conditions on the path can be relaxed. For the proof the Fredholm property and the local existence of trivializing families are essential. Thus the statement holds for a path of Dirac operators associated to an $\A$-vector bundle on a closed manifold. For a more general statement in the case $\A=\bbbc$ see \cite{rs}.  

(2) For loops related results have been proven before in geometric situations using different methods. See \cite[Theorem 3.2]{dz2}, \cite[Cor. 3]{lp2} and  \cite[Theorem 8.1]{bcp}. 

(3) For a path of families of Dirac operators this result combined with the family index theorem reproves  \cite[Theorem 0.1]{dz2} expressing the Chern character of the spectral flow in terms of $\eta$-forms of $D_0+A_0$ and $D_1+A_1$ and an integral over a local term.

{\bf Example.} 
The first part of the following discussion is a straightforward generalization of \cite[Theorem 4.4]{dz2}. 

Let $D \in RSF(H_{\A})$ and assume that there is a trivializing operator $A$ for $D$. Let $u \in \U(H_{\A})$ such that $[D,u]$ and $[D,u^{-1}]$ are bounded on $H_{\A}$, furthermore compact and densely defined as operators from $H(D)$ to $H_{\A}$. Let $P=1_{\ge 0}(D+A)$. Then $uPu^{-1}= 1_{\ge 0}(D+[u,D]u^{-1}+uAu^{-1})$ is also a generalized spectral section of $D$, hence $P-uPu^{-1}$ is compact and the Toeplitz operator $PuP:P(H_{\A}) \to P(H_{\A})$ is Fredholm with index $\ind(uPu^{-1},P)$ (see \S \ref{relindex}, Example (2)). The index does not depend on $P$ since for a second generalized spectral section $Q$ of $D$ 
\begin{eqnarray*}
\ind(uQu^{-1},Q)&=& \ind(uQu^{-1},uPu^{-1}) + \ind(uPu^{-1},P) + \ind(P,Q)\\
& =& \ind(uPu^{-1},P) \ .
\end{eqnarray*}
The index can be expressed as a spectral flow:
For $t \in [0,1]$ define $D_t:=(1-t)D+tuDu^{-1}=D+t[u,D]u^{-1}$. Then 
$$\ind(PuP+(1-P))=\ind(uPu^{-1},P)=\spfl_i((D_t)_{t \in [0,1]},A,uAu^{-1}) \ .$$
If the resolvents of $D$ are compact, then by the previous proposition this equals the index of the closure of $\ra_t +D_t + \chi_0(t)A  + \chi_1(t)uAu^{-1}$ on $L^2(\bbbr) \ten H_{\A}$. 

This example plays a role in the pairing of odd $K$-theory with $KK$-theory (see \cite{wa2}).

\section{The spectral flow and Bott periodicity}
\label{spflBott}

In this section we define the spectral flow for elements in $RSF(C([0,1],H_{\A}))$ that are invertible at the endpoints. We assign an element in $\U(C_0((0,1),K(H_{\A}))^{\sim})$ to such a path and then map it to $K_0(\A)$ via $$\beta:K_1(C_0((0,1),K(H_{\A}))) \cong K_1(C_0((0,1),\A)) \cong K_0(\A)$$ defined to be the composition of the standard isomorphism with Bott periodicity. 

Let $U:\bbbr \to S^1$ be a continuous function with $U(x)=1$ for $x \in \bbbr \setminus (-1,1)$ and such that the continuous extension $U:\ov{\bbbr} \to S^1$ has winding number one.

\begin{ddd}
Let $(D_t)_{t \in [0,1]} \in RSF(C([0,1],H_{\A}))$ and assume that $D_0$ and $D_1$ are invertible and let $\chi \in \C(\bbbr)$ be a normalizing function of $(D_t)_{t \in [0,1]}$ with $\chi(D_0)^2=\chi(D_1)^2=1$. Then $U(\chi(D_t)) \in \U(C_0((0,1),K(H_{\A}))^{\sim})$. We define
$$\spfl((D_t)_{t \in [0,1]}):=\beta[U(\chi(D_t))] \in K_0(\A) \ .$$ 
\end{ddd}

The spectral flow does not depend on the choice of $U$ and $\chi$ since different choices of $U$ and $\chi$ give rise to elements in $\U(C_0((0,1),K(H_{\A}))^{\sim})$ that can be joined by a homotopy.

The following properties are easy to verify:

(1) It is additive with respect to concatenation of paths.

(2) The restriction of $\spfl$ to $BSF(C([0,1],H_{\A}))$ is functorial in $\A$.

(3) $\spfl((D_t)_{t \in [0,1]})= \spfl((f(D_t))_{t \in [0,1]})$ for any non-decreasing function $f \in \C(\bbbr)$ with $f(0)=0$ and $f'(0)>0$.

(4) The spectral flow is homotopy invariant: If $(D_{s,t})_{s,t \in [0,1]}$ is a regular operator on $C([0,1] \times [0,1],H_{\A})$ such that $(D_{s,0})_{s \in [0,1]}$ and $(D_{s,1})_{s \in [0,1]}$ have a bounded inverse as operators on $C([0,1],H_{\A})$, then $$\spfl((D_{0,t})_{t \in [0,1]})=\spfl((D_{1,t})_{t \in [0,1]}) \ .$$

\begin{prop}
\label{spbott}
Let $(D_t)_{t \in [0,1]} \in RSF(C([0,1],H_{\A}))$ be such that $D_0,D_1$ are invertible. Then $\spfl(D_t)$ equals the image of the class of the truly unbounded odd Kasparov $(\bbbc,C_0((0,1), \A))$-module  $(H_{C_0((0,1),H_{\A})},(D_t)_{t \in [0,1]})$ under the composition
$$KK_1(\bbbc,C_0((0,1), \A)) \to K_1(C_0((0,1),\A)) \to K_0(\A) \ .$$
Here the first map is the standard isomorphism and the second is induced by Bott periodicity.
\end{prop} 

\begin{proof}
This follows from the definition of the standard isomorphism.
\end{proof}
  
Let $\chi_i \in \C(\bbbr),~i=0,1,$ be such that $\chi_0|_{(-\infty,0])}=1$ and $\chi_0|_{[\frac 13,\infty)}=0$, furthermore $\chi_1|_{[1,\infty)}=1$ and $\chi_1|_{(-\infty,\frac 23]}=0$.

\begin{lem}
\label{specflind}
Let $D \in RSF(H_{\A})$ and let $A_0$ and $A_1$ be trivializing operators for $D$. Then
$$\spfl(D + \chi_0(t) A_0+ \chi_1(t) A_1)=\ind(D,A_0,A_1) \in K_0(\A) \ .$$ 
\end{lem} 
 
\begin{proof}
Let $P=1_{\ge 0}(D+A_0)$ and $Q=1_{\ge 0}(D+A_1)$ and let $\chi$ be a normalizing function for $D$ such that $\chi(D+A_0)=2P-1$ and $\chi(D+A_1)=2Q-1$. Then $\chi(D + \chi_0(t) A_0+ \chi_1(t) A_1)-(t(2P-1)+(1-t)(2Q-1)) \in C_0((0,1),K(H_{\A}))$ by Prop. \ref{diffspec}. It follows that the class of $U(\chi(D + \chi_0(t) A_0+ \chi_1(t) A_1))$ agrees with the class of $U(t(2P-1)+(1-t)(2Q-1))$ in $K_1(C_0((0,1),\A))$. By Example (3) in \S \ref{relindex} and Prop. \ref{unique}
$$\spfl(D + \chi_0(t) A_0+ \chi_1(t) A_1) = \ind(Q,P) =\ind(D,A_0,A_1) \ .$$
\end{proof}

\begin{theorem} Let $(D_t)_{t\in [0,1]} \in RSF(C([0,1],H_{\A}))$ such that $D_0$ and $D_1$ are invertible. Assume that locally  trivializing families for $(D_t)_{t\in [0,1]}$ exist. Then
$$\spfl_i((D_t)_{t\in [0,1]},0,0)= \spfl((D_t)_{t\in [0,1]}) \ .$$
\end{theorem}

\begin{proof}
By concatenation we may restrict to the case where a trivializing family $(A_t)_{t \in [0,1]}$ for $(D_t)_{t\in [0,1]}$ exists.

We define a new path $(\tilde D_t)_{t \in [0,1]} \in RSF(C([0,1],H_{\A}))$ by setting $\tilde D_t=D_0+ \chi_1(3t)A_0$ for $t \in [0,\frac 13]$, furthermore $\tilde D_t=D_{3t-1}+A_{3t-1}$ for $t \in [\frac 13,\frac 23]$ and $\tilde D_t=D_1+\chi_0(3t-2)A_1$ for $t \in [\frac 23,1]$. Then by homotopy invariance $\spfl_i((\tilde D_t)_{t \in [0,1]} ,0,0)=\spfl_i((D_t)_{t\in [0,1]},0,0)$ and $\spfl((\tilde D_t)_{t \in [0,1]})=\spfl((D_t)_{t \in [0,1]})$. 
From the additivity of $\spfl_i$ and $\spfl$ with respect to concatenation of paths, from the fact that $\spfl_i$ and $\spfl$ vanish for paths of invertible operators and from the previous lemma it follows that
$\spfl_i((\tilde D_t)_{t \in [0,1]} ,0,0)=\spfl((\tilde D_t)_{t\in [0,1]}).$
\end{proof}

\section{Uniqueness of the spectral flow}

The justification for our definition of the spectral flow is that we use classical methods and that our notion restricts to those previously defined in less general situations. There might be other ways to obtain a generalization of the classical spectral flow. In this section we show that under some natural assumption the spectral flow is unique. 

\begin{ddd}
A spectral flow is a map $$\Spfl:\{D \in RSF(C([0,1],H_{\A}))~|~ D(0) \mbox{ and } D(1) \mbox{ have a bounded inverse} \} \to K_0(\A)$$ with the following properties:
\begin{itemize}
\item[(I)] The restriction of $\Spfl$ to $BSF(C([0,1],H_{\A}))$ is functorial in $\A$.
\item[(II)] $\Spfl(D)= \Spfl(\chi(D))$ for any normalizing function $\chi$ of $D$.
\item[(III)] It is additive with respect to direct sums of operators.
\item[(IV)] The spectral flow of a path that is symmetric with respect to the point $\frac 12$ vanishes.
\item[(V)] If $P \in K(H_{\A})$ is a projection, then $\Spfl((2tP-1)_{t \in [0,1]})=[P] \in K_0(\A)$.
\end{itemize}               
\end{ddd}

\begin{prop}
\label{axspec}
The spectral flow $\Spfl$ is uniquely defined and agrees with $\spfl$.
\end{prop}

\begin{proof} The proof proceeds in three steps.

(1) We begin by showing that the spectral flow of a loop $T \in BSF(C(S^1,H_{\A}))$ such that $T(0)$ is an involution and $T^2-1$ is compact, only depends on its class $[T] \in KK_1(\bbbc,C_0((0,1),\A)) \subset KK_1(\bbbc,C(S^1, \A))$. 

By (III) we may restrict to the case $[T]=0$.

Then there is an operator $\tilde T \in BSF(C(S^1 \times [0,1],H_{\A}\oplus H_{\A}))$ such that $\tilde T(\cdot,0)=T \oplus 1$ and $\tilde T(\cdot,1)=T(0) \oplus 1$. By functoriality the spectral flow of $T \oplus 1$ equals the spectral flow of the loop $T'$ with $T'(x)= \tilde T(0,2x)$ for $x \in [0,\frac 12]$ and $T'(x)=\tilde T(0,1-2x)$ for $x \in [\frac 12,1]$. The spectral flow of $T'$ vanishes by (IV).

(2) We show that the spectral flow is determined by its action on loops as in (1).

Let $\chi\in \C(\bbbr)$ be a normalizing function for $D$ such that $\chi(D(0))$ and $\chi(D(1))$ are involutions. Denote $P_0:= 2 \chi(D(0))-1$ and $P_1:=2\chi(D(1))-1$. By the Stabilization Theorem there are unitaries $V:P_0(H_{\A})\oplus H_{\A} \to  P_1(H_{\A}) \oplus H_{\A}$ and $W:(1-P_0)H_{\A} \oplus H_{\A} \to (1-P_1)H_{\A} \oplus H_{\A}$. We decompose $H_{\A}^3= (P_iH_{\A}\oplus H_{\A}) \oplus ((1-P_i)H_{\A}\oplus H_{\A}),~i=0,1,$ and define $\tilde P_i$ to be the projection onto $(P_iH_{\A}\oplus H_{\A})$. We define $U=V \oplus W \in {\cal U}(H_{\A}^3)$. Then $U\tilde P_0U^*=\tilde P_1$. Since $\U(H_{\A}^3)$ is contractible, we can join $U$ to the identity by a continuous path $\tilde U:[0,1] \to \U(H_{\A}^3)$ which is unique up to homotopy. Extend $\chi(D(t))$ by $1$ to $H_{\A}^3$ and define the loop $L(D) \in BSF(C(S^1,H_{\A}^3))$ by $L(D)(x)=\tilde U(2x)\chi(D_0) \tilde U(2x)^*$ for $x \in [0,\frac 12]$ and $L(D)(x)= \chi(D(2x-1))$ for $x \in [\frac 12,1]$. By construction $L(0)=L(1)$ is invertible, $L(D)^2-1 \in K(C(S^1,H_{\A}^3))$ and $\Spfl(L(D))= \Spfl(D)$.

(3) We show that $\Spfl(T)=\beta[T]$ for $T$ as in (1), where $\beta:KK_1(\bbbc,C_0((0,1), \A)) \to K_0(\A)$ is given by Bott periodicity. For that aim we evaluate condition (V). From Prop. \ref{spbott} and from $\spfl(2tP-1)=[P]$ it follows that $\beta[L(2tP-1)]=[P]$.  Hence (V) implies $\Spfl(L(2tP-1))=\beta[L(2tP-1)]$ and hence also $\Spfl(L(1-2tP))=-[P]=-\beta[L(1-2tP)]$.  By $\beta[L(2tP-1)]=[P]$ any element in $KK_1(\bbbc,C_0((0,1),\A)$ can be represented by a direct sum of loops of the form $L(2tP-1)$ resp. $L(1-2tP)$.
 
The second assertion follows from the fact that $\spfl$ fulfills the axioms. 
\end{proof}

{\bf Remark.} It is clear from the proof that the spectral flow for loops is uniquely defined by the axioms (I)-(IV) and a suitably modified axiom (V).

\section{Applications}
\label{appl}

(I) In geometric examples one often deals with families of operators with the Hilbert space also depending on the parameter. Then a rigorous definition of the spectral flow makes it necessary to choose a trivialization of the family of Hilbert spaces. We assume the Hilbert space to be separable, hence this amounts to choosing an orthonormal frame. The structure group of the bundle of orthonormal frames is (per definition) the unitary group endowed with the strong operator topology. A definition of the spectral flow that is independent of the choice of the trivialization has to allow for paths which are only strongly continuous.  In contrast to previous definitions of the spectral flow known to the author the spectral flow introduced here is invariant under conjugation by a strongly continuous family of unitary operators.
 
The following example illustrates this point:

Let $\pi:M \to S^1$ be a fiber bundle of closed manifolds endowed with a vertical Riemannian metric. Let $D$ be a vertical symmetric differential operator on $\C(M)$ such that its restriction $D_t$ to $\C(M_t)$ is elliptic for all $t \in S^1$. We get a family of unbounded Fredholm operators $D_t$ on the family of Hilbert spaces $L^2(M_t)$. There are at least two possibilities to obtain a path on a fixed Hilbert space to which we can apply the spectral flow: We define a trivialization of the bundle of vertical $L^2$-spaces by choosing a frame depending continuously on the base point. Alternatively we trivialize the fiber bundle $M$ on $S^1 \setminus \{1\}$, say, and then rescale the vertical metric on $\pi^{-1}(S^1 \setminus \{1\})$ such that the vertical volume form does not depend on the base point. Note that in the second case the transformed path of operators is not a loop any more. In the first resp. second case two different trivializations give rise to a strongly continuous loop resp. path of unitary operators. In the second case this is due to the fact that a continuous family of diffeomorphisms of a closed Riemannian manifold does in general not induce a continuous but only a strongly continuous family of operators on the $L^2$-space of functions on the manifold.
The spectral flow of both paths agrees.

(II) We discuss the spectral flow of a path of elliptic symmetric differential operators of order one acting on the sections of a hermitian vector bundle $E$ on a Riemannian manifold $M$ with boundary $N$ (compare with \cite{blp}).\\
Let $(D_t)_{t \in [0,1]}$ be a path of elliptic symmetric differential operators of order one acting on $\C(M,E)$ such that the coefficients depend in a continuous way on $t$. Assume that in a collar $U \cong [0,\ve) \times N$ the operator $D_t$ takes the form 
$$\sigma_t(x,y)(\ra_x+ E_t(x))$$ where $\sigma_t$ is a smooth skew-hermitian bundle isomorphism on $E$ with  $\sigma_t^2=-1$, and $E_t(x)$ is a differential operator on $\C(\{x\} \times N,E|_{\{x\} \times N})$ that is odd with respect to the grading on $E|_{\{x\} \times N}$ with grading operator $-i\sigma_t$. Let $P_t$ be a strongly continuous path of projections on $L^2(N,E|_N)$ with $\sigma_t P_t= (1-P_t)\sigma_t$. We assume that $P_t$ restricts to a strongly continuous path of bounded operators on $\C(N,E|_N)$. For $t \in [0,1]$ we define 
$$\dom D_t:= H^1_{P_t}(M,E):=\{f \in H^1(M,E)~|~ P_t(f|_{N})=0\} \ .$$
We assume that $D_t$ is selfadjoint and that $\C_{P_t}(M,E):=\{f \in \C(M,E)~|~ P_t(f|_{N})=0\}$ is a core for $D_t$.

Let $\phi \in \C_c(U)$ such that $\phi|_{[0,\ve/2)\times N}=1$ and such that $\phi$ is independent of $y$.

\begin{prop}
The family $(D_t)_{t \in [0,1]}$ defines a regular operator on $C([0,1],L^2(M,E))$.
\end{prop}

\begin{proof}
We extend $P_t$ to a bounded operator on $L^2(U,E|_U)$  by parallel transport in the direction orthogonal to the boundary.

If $a \in \C_{P_{t_0}}(M,E),~t_0 \in [0,1],$, then $a+ \phi(P_0-P_t)a$ is a preimage of $a$ under the evaluation map $\ev_{t_0}:\{f \in C([0,1],\C(M,E))~|~ f(t) \in \C_{P_t}(M,E) \} \to \C_{P_{t_0}}(M,E)$. It follows that the evaluation map is surjective, thus by Lemma \ref{critreg2} the operator $(D_t)_{t \in [0,1]}$ is regular on the Hilbert $C([0,1])$-module $C([0,1],L^2(M,E))$. 
\end{proof}

Assume now that $E_t$ and $\sigma_t$ do not depend on $x$ and that $P_t=1_{\ge 0}(E_t+A_t)$, where $(A_t)_{t \in [0,1]}$ is a trivializing family of $(E_t)_{t \in [0,1]}$ on $C([0,1],L^2(N,E_N))$.  The method of the proof of the following proposition is well-known to experts in family index theory.

\begin{prop}
\label{boundval}
The family $(D_t)_{t \in [0,1]}$ defines a Fredholm operator on $C([0,1],L^2(M,E))$. Furthermore it is gap continuous.
\end{prop}

\begin{proof}
We begin by constructing a parametrix $(Q_t)_{t \in [0,1]}$ of $(D_t)_{t \in [0,1]}$ on $C([0,1],L^2(M,E))$. Let $\phi_1=(1-\phi)$ and $\phi_2 =\phi$ and let $\gamma_1,\gamma_2$ be such that $\supp(\gamma_i-1)\cap \supp \phi_i=\emptyset$ and such that $\supp \gamma_1 \cap \ra M=\emptyset$ and $\supp \gamma_2 \subset U$. Let $Q^1_t$ be a local selfadjoint parametrix of $(D_t)_{t \in [0,1]}$ depending continuously on $t$ in the operator norm such that $\phi_1(Q^1_tD_t-1)\gamma_1$ and  $\gamma_1(D_tQ_t^1-1) \phi_1$ are operators with smooth integral kernels depending continuously on $t$. We define a parametrix on $U$: Let $E_N$ be the pullback of $E|_N$ to $[0,\infty) \times N$. Set $E'_t=E_t+A_t$ and consider the operator $\sigma_t(\ra_x + E'_t)$ on $L^2([0,\infty) \times N,E_N)$ with domain $$\{f \in \C_c([0,\infty) \times N,E_N) ~|~ P_t(f|_N)=0 \} \ .$$ We use the same notation for its closure. 
Let $f \in L^2([0,\infty),L^2(N,E_N))$, hence $f(x) \in L^2(N,E_N)$ for $x \in [0,\infty)$, and define 
$$(Q^2_tf)(x)=-\int_0^xe^{-(x-y)E_t'}P_t \sigma_t f(y)~dy +  \int_x^{\infty}e^{-(x-y)E_t'}(1-P_t) \sigma_t f(y)~dy \ .$$
Then $Q_t^2$ is a bounded inverse of $\sigma_t(\ra_x + E'_t)$ on $L^2([0,\infty)\times N,E_N)$ and $Q_t^2\gamma_2$ is a compact operator depending continuously on $t$. It is easily checked that $Q_t:=\phi_1Q_t^1\gamma_1+\phi_2Q_t^2\gamma_2$ is a left parametrix of $D_t$ and its adjoint is a right parametrix of $D_t$. 

Let $K_t=D_tQ_t^*-1$.
Gap continuity follows from the fact that $$(D_t+i)^{-1}=(D_t+i)^{-1}D_tQ_t^*-(D_t+i)^{-1}K_t \ .$$
The right hand side is a family of compact operators depending continuously on $t$.
\end{proof}

(IV) In the following we derive conditions for the well-definedness of the spectral flow for elliptic differential operators on noncompact manifolds. First we give an example for a path that is not gap continuous: Let $M=\bbbr$ and let $f \in C(\bbbr)$ with $f|_{[-1,1]}=1$ and $f|_{\bbbr\setminus ]-2,2[}=2$ and for $t \in [0,1]$ let $f_t \in C(\bbbr)$ be defined by $f_t(x)= f(tx)$. We define $D_t$ on $L^2(\bbbr)$ to be multiplication by $f_t$. The path $(D_t)_{t \in [0,1]}$ and its resolvents are strongly continuous but not continuous at $t=0$. \\
Now let $(D_t)_{t \in [0,1]}$ be a path of elliptic differential operators acting on the sections of a hermitian bundle $E$ on a complete Riemannian manifold $M$ with coefficients depending locally continuously on $t$. Assume  that $D_t$ is selfadjoint and that $\C_c(M,E)$ is a core for $D_t$ for each $t \in [0,1]$. Assume furthermore that there is a positive function $f \in \C_c(M)$ such that $D_t^2+f$ is invertible for all $t \in [0,1]$ with uniformly bounded inverse. It follows that $D_t$ is Fredholm for any $t$. By Lemma \ref{critreg2} the family $(D_t)_{t \in [0,1]}$ is regular on the Hilbert $C([0,1])$-module $C([0,1],L^2(M,E))$ and by Lemma \ref{critrsf} it is Fredholm since we have that  $$(D_t^2+1)^{-1}- (D_t^2+f+1)^{-1}=(D_t^2+1)^{-1}f(D_t^2+f+1)^{-1} \in C([0,1],K(L^2(M,E))$$ and $(D_t^2+f+1)^{-1}$ is uniformly bounded by a constant $c<1$. Hence the spectral flow of $(D_t)_{t \in [0,1]}$ is well-defined. 

The results in (III) and (IV) generalize to families of elliptic operators resp. elliptic operators associated to $C^*$-vector bundles.

\section{The noncommutative Maslov index of a pair of paths of Lagrangians} 
\label{lagr}

In this section we generalize the Maslov index of a pair of paths of Lagrangians to $C^*$-modules. We refer to \cite{clm} for a survey on the classical Maslov index and its geometric applications. 
Lagrangians over $C^*$-algebras were studied in \cite{wa} for example. 

Let $\A^{2n}$ be endowed with the standard $\A$-valued scalar product $$\langle x,y \rangle=\sum_{i=1}^{2n} x_i^*y_i \ .$$ Let $I=\left(\begin{array}{cc} i & 0 \\ 0 & -i \end{array}\right)$ act on $\A^{2n}=\A^n \oplus \A^n$ and define the $\A$-valued skewhermitian form $\omega(x,y)=\langle x,Iy \rangle$ on $\A^{2n}$. 

We call two submodules $L_0,L_1 \subset \A^{2n}$ {\it transverse} if $L_0 \oplus L_1=\A^{2n}$.

We call a submodule $L \subset \A^{2n}$  {\it Lagrangian} if it is isotropic with respect to $\omega$ and if $IL \oplus L = \A^{2n}$.

Lagrangian submodules are in one-to-one correspondence to unitaries in $M_n(\A)$ via the map
$$u \mapsto L(u):=\frac 12 \left(\begin{array}{cc} 1 & u^* \\ u & 1\end{array}\right)\A^{2n} \ .$$
For two Lagrangians $L(u_0),L(u_1)$ we set $$U(L(u_0),L(u_1))=\left(\begin{array}{cc} 1 & 0 \\ 0 & u_0u_1^* \end{array}\right) \in \U(\A^{2n}) \ .$$ The class $[U(L(u_0),L(u_1))]=[u_0u_1^*] \in K_1(\A)$ is invariant under the action of unitaries preserving $\omega$ on the Lagrangians. It can be characterized as follows: If $U$ is a unitary with $[U,I]=0$ and  $UL(u_1)=L(u_0)$, then $[U(L(u_0),L(u_1))]=[U]\in K_1(\A)$. The element $[U(L(u_0),L(u_1))]\in K_1(\A)$ will play a role in the definition of the difference element for Lagrangian spectral sections.

Two Lagrangian submodules $L(u_0),L(u_1)$ are transverse if and only if $u_0-u_1$ is invertible.

For $u \in \U(\A^n)$ let $D(u)$ be the operator $I \ra_x$ on $\C([0,1],\A^{2n})$ with boundary conditions $f(0)\in L(1),~f(1) \in L(u)$. The closure of $D(u)$ as an unbounded operator on $L^2([0,1],\A^{2n})$ is denoted by $D(u)$ as well. 
 
\begin{prop}
The operator $D(u)$ is regular selfadjoint with compact resolvents.
\end{prop}

\begin{proof}
If $\exp(-2 i\lambda)$ is in the resolvent set of $u$, then $(D(u) - \lambda)$ is invertible with inverse given by 
$$((D(u) - \lambda)^{-1}f)(x)=\int_0^x I e^{I\lambda(x-y)}f(y)~dy + \int_0^1
e^{I\lambda(x-y)} A(y) f(y) ~dy $$
where $$A(y)=\frac{i}{u-e^{-2i\lambda}} \left(\begin{array}{ccc}
u & e^{-2i\lambda(1-y)}\\
ue^{-2i\lambda y} & e^{-2i \lambda}
\end{array}\right) \ .$$
\end{proof}

We need the following technical lemma:

\begin{lem}
If $L(u(t))$ is a continuous path of Lagrangians, then
there is  $W \in \C([0,\frac 13],C([0,1],\U(\A^{2n}))$ with $[W,I]=0$, furthermore $W(0,t)L(u(t))=L(u(0))$, $W(x,0)=1$ and $W(\frac 13,t)=1$ for all $x,t$. 
\end{lem}

\begin{proof}
We give an algorithm for the construction of $W$: Choose $c=0$ if possible, else $c>1$, such that $\|u(t)-u(ct)\| \le \frac 12$ for all $t \in [0,1]$. Then the spectrum of $U(L(u(ct),L(u(t))$ in $M_{2n}(C([0,1],\A))$ is not equal to $S^1$. Hence there is a path $W \in \C([0,\frac 13],U(C([0,1],\A)^{2n})$ with $W(0,t)=U(L(u(ct),L(u(t))$ and $W_1(\frac 13,t)=1$. Now repeat the procedure for $u(ct)$. After finitely many steps $c$ may be chosen $0$.  
\end{proof}

Let $(L(u_0(t)),L(u_1(t)))_{t \in [0,1]}$ be a pair of continuous paths of Lagrangian subspaces of $\A^{2n}$ and assume that $L(u_0(i)),~L(u_1(i))$ are transverse for $i=0,1$. 

Consider the operator $I \ra_x$ on $\C([0,1],\A^{2n})$ with boundary conditions $f(0)\in L(u_0(t)),~f(1) \in L(u_1(t))$. Its closure as an unbounded operator on $L^2([0,1],\A^{2n})$ is denoted by $D_t$. Then $(D_t)_{t \in [0,1]} \in RSF(L^2([0,1],C([0,1],\A)^{2n}))$, where we consider $L^2([0,1],C([0,1],\A)^{2n})$ as a Hilbert $C([0,1],\A)$-module. Furthermore $(D_t)_{t \in [0,1]}$ has compact resolvents, and $D_0, D_1$ are invertible. 

We define the Maslov index of $(L(u_0(t)),L(u_1(t)))_{t \in [0,1]}$ by
$$\mu(L(u_0),L(u_1)):= \spfl((D_t)_{t \in [0,1]}) \in K_0(\A) \ .$$

The results of previous sections provide alternative ways of describing the Maslov index:

By the previous lemma we find $W_t \in \C([0,1],\U(\A^{2n}))$ depending continuously on $t$ with $W_t(i)L(u_i(t))=L(u_i(0)),~i=0,1,$ and $[W_t,I]=0$. Then 
$$W_t D_t W^*_t=D_0 +  W_tI(\ra_x W^*_t) $$
is a family of operators with fixed domain $\dom D_0$.  

By Lemma \ref{boundspec} there are paths of generalized spectral sections for $W_tD_tW_t^*$, hence
 $$\mu(L(u_0),L(u_1))= \spfl_i(W_t D_t W^*_t,0,0) \ .$$

We extend a path $u:[0,1] \to \U(\A^n)$ to $\bbbr$ by setting $u(t)=u(0)$ for $t<0$ and $u(t)=u(1)$ for $t>1$. We call $u$ resp. $L(u)$ regular if the extended path is $C^1$.

If $u_0,u_1:[0,1] \to {\cal U}(\A^n)$ are regular and $u_0(i)-u_1(i),~i=0,1$ invertible, then we define ${\cal D}(L(u_0),L(u_1))$ to be the closure of the operator $\ra_t+I\ra_x$ on $L^2(\bbbr \times [0,1],\A^{2n})$ with domain $$\{f\in \C_c(\bbbr \times [0,1] ,\A^{2n})~|~f(t,0)\in L(u_0(t)),~f(t,1) \in L(u_1(t)) \}  \ .$$ 
By Prop. \ref{spflind}
$$\mu(L(u_0),L(u_1))=\ind ~{\cal D}(L(u_0),L(u_1)) \ .$$

The following properties are elementary to prove:

We have that $$\mu(L(u_0),L(u_1))=-\mu(L(u_1),L(u_0)) \ .$$

Furthermore the Maslov index is functorial in $\A$: If $f:\A \to \B$ is a unital $C^*$-algebra homomorphism, then $$\mu(L(f(u_0)),L(f(u_1)))=f_*\mu(L(u_0),L(u_1)) \in K_0(\B) \ .$$

The Maslov index is homotopy invariant:
If $u_i(x,t) \in C([0,1]\times [0,1],\U(\A^n)),~i=0,1,$ with $u_0(x,i)-u_1(x,i)$ invertible for $i=0,1$, then $$\mu(L(u_0(0, \cdot)),L(u_1(0,\cdot))=\mu(L(u_0(1, \cdot)),L(u_1(1,\cdot))$$
by the functoriality of the Maslov index and by $\ev_0=\ev_1:K_0(C([0,1],\A))\to K_0(\A)$.  

\subsection{The Maslov index and the triple Maslov index}
The noncommutative Maslov index of a triple of pairwise transverse Lagrangians $(L_0,L_1,L_2)$ is a generalization of the classical triple Maslov index. It was introduced in \cite{bk} for families and in \cite{wa} for $C^*$-algebras. 

We recall its definition:

Denote by $P_i \in M_{2n}(\A)$ the orthogonal projection onto $L_i$.
 
The hermitian $\A$-valued form
$$h:L_0 \times L_0 \to \A,~ (v,w) \mapsto \langle v_2,w_1\rangle \ ,$$ where $v=v_1+v_2,~w=w_1+w_2$ is the decomposition 
with respect to $\A^{2n}=L_1\oplus L_2$, is nonsingular since the associated selfadjoint matrix is given by 
$$P_0(P_1+P_2)^{-1}P_2IP_1(P_1+P_2)^{-1}P_0$$
$$=(P_0+P_1)(P_1+P_2)^{-1}(P_1+P_2)IP_1(P_1+P_2)^{-1}(P_0+P_2) \ ,$$ which has closed range.

The class of $h$ in $K_0(\A)$ is denoted by $\tau(L_0,L_1,L_2)$ and called the Maslov (triple) index of $(L_0,L_1,L_2)$.

The remainder of this section is devoted to the proof of the following proposition, which generalizes a classical formula relating the Maslov triple index to the Maslov index for pairs. 

\begin{prop}
If $(L_0(t),L_1(t),L_2(t))_{t \in [0,1]}$ is a triple of paths of Lagrangians that are pairwise transverse at $t=0$ and $t=1$, then 
 $$\tau(L_0(1),L_1(1),L_2(1)) - \tau(L_0(0),L_1(0),L_2(0)) =$$
 $$2\bigl(\mu(L_0,L_1)+ \mu(L_1,L_2) + \mu(L_2,L_0)\bigr) \ .$$
\end{prop}

Without loss of generality we may assume that the paths are regular.

The proof of the formula is a relative $K$-theoretic index theorem. It uses that $\mu(L_i,L_j)$ is the index of the operator ${\cal D}(L_i,L_j)$ on $L^2(\bbbr\times [0,1],\A^{2n})$ as shown in the previous section and that also the Maslov triple index can be expressed as the index of a Dirac operator on a manifold with boundary with local boundary conditions defined by the Lagrangians. Such an operator was constructed for the  Maslov triple index for families in \cite{bk}. The (straightforward) generalization to the noncommutative context can be found in \cite{wa}.

Giving only the details needed in the proof of the formula we outline the construction of the operator whose index is the Maslov index $\tau(L_0,L_1,L_2)$.

\begin{figure}
\includegraphics{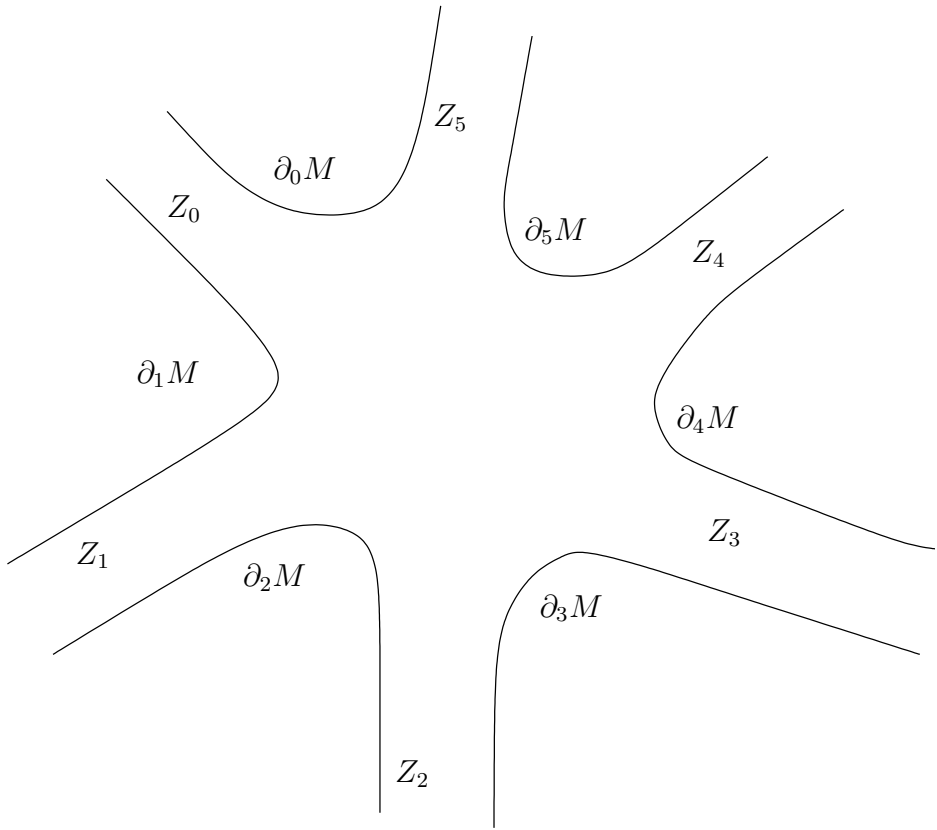}
\caption{The manifold $M$}
\end{figure}

The underlying manifold is a simply connected oriented two dimensional spin manifold $M$ with six boundary components $\ra_k M,~k=0,\dots ,5,$ isometric to $\bbbr$ and six cylindric ends $Z_k,~k=0, \dots, 5,$ isometric to $\bbbr^+ \times [0,1]$. The diffeomorphism structure of the manifold is illustrated by figure 1, which also serves to fix notation. 
We assume the metric to be of product type in a neighborhood of the boundary. Let $S$ be the spinor bundle on $M$ and $E=S \ten ((\bbbc^+)^n \oplus (\bbbc^-)^n)$ with the obvious Dirac bundle structure. Clifford multiplication with the volume form on $M$ induces a skewhermitian form on $E$. There is a connected flat region $F$ containing the boundary components and the cylindric ends and there is an orthonormal frame of  $S|_F$ invariant under parallel transport. We use it to identify $E^+|_{Z_k}$ with the trivial Dirac bundle $Z_k \times \bbbc^{2n}$ and $E^+|_{\ra_k M}$ with the trivial Dirac bundle $\ra_k M \times \bbbc^{2n}$ such that Clifford multiplication with the volume form is given by the matrix $I$.

Then the closure of $\dira^+:\dom \dira^+ \to L^2(M,E^-\ten \A)$ with domain 
$$\dom \dira^+=\{f \in \C_c(M,E^+\ten \A)~|~f(x)\in L_i \mbox{ if } x \in \ra_iM\cup \ra_{i+3}M,~i=0,1,2 \}$$
is Fredholm with index $\tau(L_0,L_1,L_2)$.

For the proof of the formula we construct two boundary value problems, namely Dirac operators with boundary conditions which will be denoted by $D_1,D_2$. 

The underlying manifold of $D_1$ is denoted by ${\cal M}_1$ and is the disjoint union of two copies $M_0,M_1$ of $M$ ($M_0$ for the terms  $\tau(L_0(0),L_1(0),L_2(0)$ and $M_1$ for $\tau(L_0(1),L_1(1),L_2(1))$ in the above formula) and six cylinders (one for each term on the right hand side). The cylindric ends of $M_0,M_1$ are denoted by $Z_i^0,Z_i^1,~i=0,\dots, 5$. On $M_1$ we reverse the orientation. On $M_0$ resp. $M_1$ we define $D_1^+=\dira^+$ as above with boundary conditions associated to the triple $(L_0(0),L_1(0),L_2(0))$ resp. $(L_0(1),L_1(1),L_2(1))$. To each term $\mu(L_i,L_j)$ in the formula we associate exactly one cylinder, on which we set $D_1^+={\cal D}(L_i,L_j)$. As usual let $D_1^-$ be the adjoint of $D_1^+$ and $D_1=D_1^+ \oplus D_1^-$.

The underlying manifold ${\cal M}_2$ of $D_2$ is the disjoint union of twelve cylinders and a compact manifold $M'$ with boundary. The manifold $M'$ is constructed as follows: We cut the cylindric ends $Z_i^0,Z_i^1,~i=0,\dots, 5$ of $M_0,M_1$ such that we obtain two incomplete manifolds $M^c_0,M^c_1$ with truncated ends $Q_i^0,Q_i^1,~i=0, \dots 5,$ isometric to $]0,1[ \times [0,1]$. Then we glue $M^c_0$ to $M^c_1$ via orientation preserving isometries $Q^0_i \cong Q_i^1$. The resulting manifold is $M'$. From the Dirac operators associated to $E \ten \A$ on $M_0,M_1$ we obtain a Dirac operator on $M'$. The manifold $M'$ has six boundary components diffeomorphic to $S^1$. We use the boundary condition on $M_j,~j=0,1$ from above to define boundary conditions on $M_j - \cup_i Q^j_i$. For $i=1,\dots 5$ we define the boundary conditions on the image of $Q_i^0$ in ${\cal M}_2$ using the path $(L_{i\Mod 3}(t),L_{(i+1)\Mod 3}(t))$ such that the resulting boundary conditions are $C^1$. The Dirac operator with these boundary conditions is the restriction of the operator $D_2^+$ to $M'$.

Furthermore we assume that the restriction of $D_2^+$ to exactly two of the twelve cylinders equals ${\cal D}(L_{i}(0),L_{i+1}(0))$ resp. ${\cal D}(L_{i}(1),L_{i+1}(1)),~ i \in \bbbz/3$. 

Again $D_2^-$ is the adjoint of $D_2^+$ and $D_2=D_2^+ \oplus D_2^-$.

By cutting and pasting along suitable cross-sections of the cylindric ends $D_1$ can be transformed into $D_2$. The $K$-theoretic relative index theorem proven in \cite{bu} for manifolds without boundary extends without difficulty to include local boundary conditions showing that $$\ind D^+_1= \ind D^+_2 \ .$$ 
Furthermore by construction
\begin{eqnarray*}
\ind  D^+_1&=&  2(\mu(L_0,L_1)+ \mu(L_1,L_2) + \mu(L_2,L_0))\\
&&- \tau(L_0(1),L_1(1),L_2(1)) + \tau(L_0(0),L_1(0),L_2(0))  \ .
\end{eqnarray*}

It remains to show that $\ind D_2^+=0$.

If the paths $L_0,L_1,L_2$ are constant, then the right hand side of the previous equation vanishes, hence in that case the formula implies that $\ind D^+_1=0$, thus also $\ind D^+_2=0$. In general $\ind D^+_2=0$ since the boundary conditions on $M'$  can be deformed into boundary conditions coming from constant paths and the contribution from the cylindric components of ${\cal M}_2$ vanishes.

\subsection{The Maslov index and Bott periodicity}
\label{maslbottper}

Let $x \in \U(\bbbc^n)$ such that $1-x$ is invertible and let $L_0:=L(x)$.

The Maslov index induces a map   $$\mu:\U((C_0((0,1),\A)^{\sim})^n) \to  K_0(\A),~u \mapsto \mu(L_0,L(u)) \ .$$
It does not depend on the choice of $L_0$ since $L_0$ can always be deformed into $L(-1)$ via a path of Lagrangians transverse to $L(1)$.
 
The map is compatible with the inductive limit with respect to $n$ and homotopy invariant. Hence we get an induced map
$$\mu:K_1(C_0((0,1), \A)) \to K_0(\A) \ .$$

\begin{prop}
The map 
$$\mu:K_1(C_0((0,1),\A)) \to K_0(\A)$$ equals the Bott periodicity map.
\end{prop}

The proposition will follow from Prop. \ref{standiso} below.

\begin{cor}
\label{maslbott}
Let $u_i \in \U(C(S^1,\A^{2n})),~i=0,1,$ be such that $L(u_0(1))$ and $L(u_1(1))$ are transverse. Then 
$$\mu(L(u_0),L(u_1))= \beta [U(L(u_1),L(u_0))]=\beta[u_1u_0^*] \in K_0(\A) .$$
\end{cor}

For the proof of the proposition note the following:

By definition of $\mu$ and by Prop. \ref{spbott} the element $[L^2([0,1],C(S^1, \A)^{2n}),D(u)] \in KK_1(\bbbc,C_0((0,1),\A))$ is mapped 
  to $\mu(L_0,L(u)) \in K_0(\A)$ under the Bott periodicity map. 
 
Hence it remains to show that the image of $[L^2([0,1],C(S^1, \A)^{2n}),D(u)]$ is $[u] \in K_1(C_0((0,1), \A))$ under the standard isomorphism  $KK_1(\bbbc,C_0((0,1))\ten \A) \cong K_1(C_0((0,1))\ten \A)$. In the remainder of this section we will show the more general fact that for any unitary $u \in \U(\A^n)$ the element $[D(u)] \in KK_1(\bbbc,\A)$ corresponds to $[u] \in K_1(\A)$ under the standard isomorphism. 

For $u \in \U(\A^n)$ let  $D'(u)$ be the closure of the operator $i \ra_x$ on $L^2([0,1],\A^n)$ with domain $\{f \in \C([0,1],\A^n)~|~f(0)=uf(1)\}$. The operator $D'(u)$ is the Dirac operator $i \ra_x$ on $S^1$ twisted by the bundle $V(u):=[0,1] \times \A^n /(0, v) \sim (1,u^* v)$ endowed with the trivial connection. 

\begin{lem}
The element $[D'(1)] \in KK_1(C_0((0,1)), \bbbc)$ corresponds to $1 \in KK_0(\bbbc,\bbbc)$ under Bott periodicity. Hence the Kasparov product $\ten [D'(1)]$ induces the  isomorphisms
$K_0(C_0((0,1), \A)) \cong K_1(\A)$ and $K_1(C_0((0,1),\A)) \cong K_0(\A)$.
\end{lem}

\begin{proof}
The assertion follows the fact that
the Kasparov product
$$\otimes [D'(1)]: K_1(C_0((0,1))) \to K_0(\bbbc)$$
maps the generator $[u(t)=e^{2 \pi i t}]$ of $K_1(C_0((0,1)))$ to $$\ind \bigr[1_{\ge 0}(D'(1))u 1_{\ge 0}(D'(1)) +(1-  1_{\ge 0}(D'(1)))\bigl]=1 \ .$$  
\end{proof} 

\begin{prop}
The map $K_1(\A) \to   KK_1(\bbbc, \A),~[u] \mapsto [D'(u)]$ is the standard isomorphism.
\end{prop}

\begin{proof}
It is enough to show that the map equals the composition
$$K_1(\A) \to K_0(C(S^1,\A)) \to KK_0(\bbbc,C(S^1,\A)) \stackrel{\ten [D'(1)]}{\longrightarrow} KK_1(\bbbc,\A) \ .$$
Here we define the first map by 
$[u] \mapsto [V(u)]$. Note that for $u \in \U(\A^n)$ the map $[u] \mapsto [V(u)]-[S^1 \times \A^n] \in K_0(C_0((0,1),\A))$ is the standard isomorphism $K_1(\A) \to K_0(C_0((0,1),\A))$. The second map is defined to be the standard isomorphism. 

Let $c$ be an odd selfadjoint generator of the Clifford algebra $C_1$ with $c^2=1$. The element $[D'(u)]$  corresponds to $[cD'(u)]:=[C_1 \ten L^2([0,1],\A),cD'(u)]$ under $KK_1(C(S^1),\A) \cong KK_0(C(S^1),C_1\ten \A)$.
The Kasparov product $[V(u)]\ten [D'(1)]$ is defined as the image of the Kasparov product $[V(u)]\ten [cD'(1)]$ under  $KK_0(\bbbc, C_1\ten \A) \cong KK_1(\bbbc,\A)$. The Kasparov product $[V(u)]\ten [cD'(1)]$ equals $[cD'(u)]$. We obtain that the above composition maps $[u]$ to $[D'(u)]$.   
\end{proof}

\begin{prop} 
\label{standiso}
The operators $D(u)$ and $D'(u)$ define the same class in $KK_1(\bbbc,\A)$, hence the map $K_1(\A) \mapsto  KK_1(\bbbc, \A),~[u] \mapsto [D(u)]$ is the standard isomorphism.
\end{prop}

\begin{proof}
We prove that in $KK_1(\bbbc,\A)$ 
$$[D(u)]+ [D'(1_n)\oplus (-D'(1_n))]- [D'(1_n)\oplus (-D'(u^*))]-[D(-1_n)]=0 \ .$$
Then the first assertion follows from $[D'(1_n)]=[D(-1_n)]=0$ and $[-D'(u^*)]=[D'(u)]$ and the second from the previous proposition.

The equation is a relative $K$-theoretic index theorem. 
Let 
$${\cal G}=\diag(D(u),D'(1_n)\oplus (-D'(1_n)),D'(1_n)\oplus (-D'(u^*)),D(-1_n))$$ on ${\cal H}=L^2([0,1],(\A^{\sim})^{2n})^2 \oplus L^2([0,1],(\A^-)^{2n})^2$. Let $U:=\diag(1,u)$ and let $\chi_1, \chi_2, \chi_3 \in \C([0,1])$ with $\supp \chi_i \in [\frac{i-1}{4},\frac{i+1}{4}]$ and with $\chi_1^2+\chi_2^2+\chi_3^2=1$. Define
$$X=\left(\begin{array}{cccc} 0 & 0 &-\chi_2  &-\chi_1-U\chi_3   \\
0 & 0 & -\chi_1-U^* \chi_3  & \chi_2 \\
-\chi_2 & -\chi_1 - U \chi_3 & 0 & 0 \\
-\chi_1 - U^* \chi_3 &  \chi_2 & 0 & 0 
\end{array}\right) \ .$$

Then $X^2=1,~X=X^*$ and $X$ is odd. Furthermore $[X,{\cal G}]$ is bounded. Now the argument is as in the proof of Lemma \ref{relind} (5): We have that $[{\cal H},{\cal G}] \in \Ima(i^*:KK_1(C_1,\A) \to KK_1(\bbbc,\A))=\{0\}$.
\end{proof}

For later use we note still another possibility of representing a class in $KK_1(\bbbc,\A)$:

\begin{cor}
\label{standiso2}

Let $L(u_0),L(u_1) \subset \A^{2n}$ be two Lagrangian submodules and let $P_0,P_1$ be orthogonal projections onto them. Then the class of the closure $D(P_0,P_1)$ of
$I(\ra_t + \chi_0(t) (2P_0-1)+\chi_1(t)(2P_1-1))$ with domain ${\cal S}(\bbbr,\A^{2n})$ on $L^2(\bbbr,\A^{2n})$ in $KK_1(\bbbc,\A)$ equals $[D(u_1u_0^*)]$.
\end{cor}

\begin{proof}
We may assume that $u_0=1_n$ and set $u=u_1$.
 
Note that the closure of $\ra_t+(2P_i-1),~i=0,1$ is invertible, hence $D(P_0,P_1) \in RSF(L^2(\bbbr,\A^{2n}))$. Furthermore if $P_0,P_1$ are orthogonal, then $[D(P_0,P_1)]$ is in the image of $KK_1(\bbbc,\bbbc)\to KK_1(\bbbc,\A)$, hence it vanishes. Analogously to the previous proof one shows that
$$[D(P_0,P_1)]+ [D(1_n)\oplus (-D(1_n))]- [D(1_n)\oplus (-D(u^*))]-[D(P_0,1-P_0)]=0 \ .$$
\end{proof}

{\bf Remarks.} 
(1) The corollary shows that  elements of $K_1(\A)$ can be represented by pairs of involutions -- in this case  $((2P_0-1),(2P_1-1))$) -- which are odd with respect to the grading defined by the operator $-iI$. This is a special case of a description of $K$-theory due to Karoubi \cite{k}.

(2) The following situation will be relevant for the definition of the odd spectral flow in \S \ref{gradspec}: Let ${\cal P}$ be a $\bbbz/2$-graded projective $\A$-module endowed with an $\A$-valued scalar product. Let $\sigma$ be the grading operator. The skewadjoint operator $I=i\sigma$ induces skewhermitian form on ${\cal P}$. We call an isotropic submodule $L \subset {\cal P}$ Lagrangian if $L \oplus IL= {\cal P}$. Lagrangian submodules exist if ${\cal P^+}\cong {\cal P^-}$. Every Lagrangian submodule is of the form $L(u):=\Ran \frac 12\left(\begin{array}{cc} 1 & u^* \\ u & 1 \end{array}\right)$ for $u \in {\cal U}({\cal P}^+,{\cal P^-})$. Define $U(L(u_0),L(u_1))=\left(\begin{array}{cc} 1 & 0 \\ 0 & u_0u_1^* \end{array}\right) \in {\cal U}^+({\cal P})$. The class $[U(L(u_0),L(u_1))] \in K_1(\A)$ is invariant under the action of ${\cal U}^+({\cal P})$ on the Lagrangians.\\
The definitions and results of the previous two sections can be generalized to this situation.

\section{The odd spectral flow}

\label{gradspec}
\subsection{Generalized spectral sections}

Let $\sigma$ be the grading operator on $\hat H_{\A}=H_{\A}^+ \oplus H_{\A}^-$.
Note that $I=i\sigma$ induces an $\A$-valued resp. $B(H_{\A})$-valued skewhermitian form on $\hat H_{\A}$ resp. on $B(H_{\A}^+) \oplus B(H_{\A}^-)$.

We call a projection $P \in B(\hat H_{\A})$ Lagrangian if the involution $2P-1$ is odd. This is equivalent to $P(B(H_{\A}^+) \oplus B(H_{\A}^-))\subset B(H_{\A}^+) \oplus B(H_{\A}^-)$ being Lagrangian in the sense of \S\ref{lagr}.  

We call a spectral section that is a Lagrangian projection a Lagrangian spectral section. 

Let $D$ be an odd regular selfadjoint operator on $\hat H_{\A}$ with compact resolvents and with $[D]=0$ in $KK_0(\bbbc,\A)$. If there exists a spectral section of $D$, then there exists also a Lagrangian spectral sections of $D$ (see below for a proof).
 
Let $\chi$ be a spectral cut with $\chi(x)=0$ for $x \le 0$. It is easy to show that $\Ran \chi(D)$ is isotropic, hence the range of any spectral section $P$ with $P\chi(D)=P$ is isotropic. Furthermore $1-\sigma P \sigma$ is a spectral section as well.

If $D$ admits spectral sections, then  we define the difference element of a pair $(P_1,P_2)$ of Lagrangian spectral sections  of $D$ as follows: Let $R$ be a spectral section with $RP_1=R$ and $RP_2=R$. Then $P_1-R$ and $P_2-R$ are projections onto Lagrangian subspaces in 
$\Ran(1-\sigma R \sigma-R)$. We write $U(P_1-R,P_2-R)$ for $U(\Ran(P_1-R),\Ran(P_2-R))$ as defined in the remark at the end of \S \ref{maslbottper}. The class $[U(P_1-R,P_2-R)]\in K_1(\A)$ does not depend on $R$, hence we may set
$$[P_1-P_2]:=[U(P_1-R,P_2-R)] \in K_1(\A) \ .$$
The definition of the difference element of Lagrangian spectral sections in \cite{lp2} is different. It will follow from the splitting theorem in \S \ref{split} that both definition coincide.

\begin{ddd}
A {\rm relative odd index} of projections is a map $\ind^o$ assigning to a pair $(P_1,P_2)$ of Lagrangian projections on $\hat H_{\A}$ with $P_1-P_2 \in K(\hat H_{\A})$ an element in $K_1(\A)$ and fulfilling the following axioms: 
\begin{itemize}
\item[(I)] (Additivity.)  
$$\ind^o(P_1,P_3)=\ind^o(P_1,P_2)+\ind^o(P_2,P_3) \ .$$
\item[(II)] (Functoriality.) If $\phi:\A \to \B$ is a unital $C^*$-homomorphism, then $$\phi_*\ind^o(P_1,P_2)=\ind^o(\phi_*P_1,\phi_*P_2) \ .$$ 
\item[(III)] (Stabilization.) If $R$ is a Lagrangian projection of $\hat H_{\A}$, then $$\ind^o(P_1,P_2)=\ind^o(P_1\oplus R,P_2\oplus R)$$ where for the definition of the right hand side we choose an even isomorphism $\hat H_{\A}^2 \cong \hat H_{\A}$. 
\item[(IV)](Normalization.) If $P_1,P_2$ are Lagrangian projections and there is a projection $R$ with $RP_i=R,~i=1,2,$ and $R-P_i \in K(\hat H_{\A})$, then $\ind^o(P_1,P_2)=[U(P_1-R,P_2-R)]$. 
\end{itemize}
\end{ddd}

Any odd relative index is homotopy invariant for strongly continuous paths of Lagrangian projections and agrees for Lagrangian spectral sections of an operator admitting spectral sections with the difference element.

In the following $\beta:KK_0(\bbbc,C_0((0,1),\A)) \to K_1(\A)$ is the map induced by Bott periodicity and the standard isomorphism.

{\bf Examples:}

Let $P,Q \in B(\hat H_{\A})$ be Lagrangian projections with $P-Q$ compact. \\
(1) We identify $H:=P(\hat H_{\A})$ with $(1-P)\hat H_{\A}$ via $\sigma$ and set $\B=K(H)^{\sim}$. Then $i \sigma=\left(\begin{array}{cc}0&i\\ i&0\end{array}\right)$ defines a skewhermitian form on $\B^2$ and $P,Q \in M_2(\B)$ are projections onto Lagrangian submodule of $\B^2$.  Set $\ind_1^o(P,Q)=[U(P(\B^2),Q(\B^2))] \in K_1(\B) \cong K_1(\A)$. Functoriality and normalization are clear, additivity follows from $[U(P(\B^2),Q(\B^2))][U(Q(\B^2),R(\B^2))]=[U(P(\B^2),R(\B^2))]$.   
    
(2) Define $\ind_2^o(P,Q)$ as the image of $[F_t:=(1-t)(2Q-1) + t(2P-1)] \in KK_0(\bbbc,C_0((0,1),\A))$ under $\beta$. Property (I) holds since in $KK_0(\bbbc,C_0((0,1),\A))$ concatenation of paths corresponds to addition. Let now $P,Q,R$ as in (IV) and denote ${\cal P}=\Ran (1-\sigma R \sigma - R)$. The class $[F_t] \in KK_0(\bbbc,C_0((0,1),\A))$ equals the class of the restriction of $F_t$ to $C_0((0,1),{\cal P})$. Let $U:=U(P-R,Q-R)$ and let ${\cal V}(U^*)=\{f \in C([0,1],{\cal P}),~f(0)=U^*f(1) \}$. The image of $[F_t]$ under the inclusion $KK_0(\bbbc,C_0((0,1),\A)) \to KK_0(\bbbc,C(S^1,\A))$ is represented by the operator $(1-t)(2Q-1) + t(2P-1)$ on the $\bbbz/2$-graded $C(S^1,\A)$-module ${\cal V}(U^*)$. Since ${\cal V}(U^*)$ is a projective $C(S^1,\A)$-module, the image of $[F_t]$ in $K_0(C_0((0,1),\A))$ is $[{\cal V}(U^*)^+]-[{\cal V}(U^*)^-]$, which corresponds to $[U]\in K_1(\A)$ under Bott periodicity.
\\

In the following we introduce the odd spectral flow along the lines of the discussion of the spectral flow. Since many arguments are analogous, we focus on the modifications.

Let $D \in RSF^-(\hat H_{\A})$. 

We assume trivializing operators and families to be odd. 

If $A$ is a trivializing operator for $D$, then the generalized spectral section $1_{\ge 0}(D+A)$ is Lagrangian. If $\ind^o$ is a relative odd index and $A_0,A_1$ are two trivializing operators of $D$, then by Prop. \ref{diffspec} $$\ind^o(D,A_0,A_1):=\ind^o(1_{\ge 0}(D+A_1),1_{\ge 0}(D+A_0))$$ is well-defined.

Let $A$ be a trivializing operator for $D$ and let $\phi_R:\bbbr \to [0,1]$ be a smooth even function with $\phi_R(x)=1$ for $|x|\le R$ and $\phi_R(x)=0$ for $|x|>R+1$. For $R$ big enough $\phi_R(D)A\phi_R(D)$ is also a trivializing  operator for $D$. If $D$ has compact resolvents, then the associated generalized spectral section is a Lagrangian spectral section of $D$.

Let $(D_t)_{t \in [0,1]} \in RSF^-(C([0,1],\hat H_{\A}))$ be a path for which locally trivializing families exist and let $A_0,A_1$ be trivializing operators for $D_0,D_1$.
The odd spectral flow $\spfl_i^o((D_t)_{t \in [0,1]},A_0,A_1)$ is defined by replacing $\ind$ in the definition of $\spfl_i$ in \S \ref{spflone} by $\ind^o$. Analogues of the results of \S \ref{spflone} hold.

In order to study the relation of the odd spectral flow and the index we extend a family $(D_t)_{t \in [0,1]}$ to $\bbbr$ by setting $D_t=D_0$ for $t \in (-\infty,0]$ and $D_t=D_1$ for $t \in [1,\infty)$.

Assume now that $(D_t)_{t \in [0,1]} \in RSF^-(C([0,1],\hat H_{\A})$ has compact resolvents and that $$\bbbr \to B^-(\hat H_{\A}),~ t \mapsto (D_t-D_0)$$ is well-defined and strongly continuous, furthermore that $D_0,D_1$ are invertible. As in the proof of Lemma \ref{relind} one shows:

The closure ${\cal D}$ of the operator
$i\sigma (\ra_t +  D_t)$ with domain $\C_c(\bbbr, \dom D_0)$ is in $RSF(L^2(\bbbr,\hat H_{\A}))$. If $D_t$ is invertible for all $t \in \bbbr$, then $\ind^o({\cal D})=0$. (Recall that  $\ind^o({\cal D})$ denotes the image of $[{\cal D}] \in KK_1(\bbbc,\A)$ under the standard isomorphism $KK_1(\bbbc,\A) \cong K_1(\A)$.) The index of a concatenation is the sum of the indices of the parts. 

In the following $i\sigma(\ra_t +D_t)$ means the closure.

\begin{prop}
\label{indspecsecodd}
Let $D\in RSF^-(\hat H_{\A})$ have compact resolvents
and let $A_0,A_1$ be trivializing operators of $D$. 
Let $\chi_i \in \C(\bbbr),~i=0,1,$ with $\chi_0|_{(-\infty,-1])}=1$ and $\chi_0|_{[0,\infty)}=0$, furthermore $\chi_1|_{[2,\infty)}=1$ and $\chi_1|_{(-\infty,1]}=0$.
Then
$$\ind^o(i\sigma(\ra_t + D+\chi_0(t) A_0 + \chi_1(t) A_1))=\ind^o(D,A_0,A_1) \in K_1(\A) \ .$$  
\end{prop}

\begin{proof} We first assume that $D$ admits spectral sections.

We may assume that $P_i=1_{\ge 0}(D+A)$ is a Lagrangian spectral section. Let $Q$ be a spectral section of $D$ with $QP_i=Q$ and such that $QDQ$ is invertible on $\Ran Q$. Let $R:=\sigma(1-Q)\sigma$. Then $R$ is a spectral section of $D$ and $(1-R)D(1-R)$ is invertible on $\Ran(1-R)$. Furthermore $R P_i=P_i$. Assume that $$D+A_i=QDQ \oplus (P_i-Q) \oplus (P_i-R) \oplus (1-R)D(1-R) \ .$$
Then
\begin{eqnarray*}
\lefteqn{\ind^o(i\sigma(\ra_t + D + \chi_0(t) A_0+ \chi_1(t) A_1)}\\
&=& \ind^o\bigl(i\sigma(\ra_t + \chi_0(t) (2P_0-1)+\chi_1(t)(2P_1-1)|_{L^2(\bbbr,(R-Q)\hat H_{\A})}\bigr) \ .
\end{eqnarray*} 
Now the assertion follows from Cor. \ref{standiso2}.

If $D$ does not admit spectral sections, then define the odd operator $$\tilde d=\left(\begin{array}{cc} 0 & d \\ d & 0 \end{array}\right):\hat H_{\A} \to \hat H_{\A} \ $$ with $d$ as in the end of \S \ref{specsec}. Then $\tilde d \oplus D$ admits spectral sections and $\tilde P_i:=1_{\ge 0}((D+A_i) \oplus \tilde d)$ is a Lagrangian spectral section. By the first part
$$\ind^o(i\sigma(\ra_t +\tilde d\oplus (D + \chi_0(t) A_0+ \chi_1(t) A_1))=\ind^o(\tilde P_1,\tilde P_0) \ .$$
The assertion follows from
$\ind^o(i\sigma(\ra_t + \tilde d))=0$ and $\ind^o(\tilde P_1,\tilde P_0)=\ind^o(P_1,P_0)$.  
\end{proof}

From the previous proposition we conclude in analogy to Prop. \ref{spflind}
that if $(D_t)_{t \in [0,1]}$ is a path in $RSF^-(\hat H_{\A})$ such that $D_t$ has compact resolvents and $D_t-D_0$ is bounded and strongly continuous in $t$ and $A_0,A_1$ are trivializing operators for $D_0,D_1$, then  
$$\spfl^o_i((D_t)_{t \in [0,1]},A_0,A_1)=\ind^o(i\sigma(\ra_t +D_t + \chi_0(t) A_0 + \chi_1(t) A_1)) \in K_1(\A) \ .$$

The conditions on the path can be relaxed in concrete situations.
By using the odd family index theorem \cite{mp2} this result allows one to express the odd spectral flow of a path of families of Dirac operators in terms of $\eta$-forms of the endpoints and an integral over a local term. \\

\begin{prop}
\label{uniqueodd}
The relative odd index $\ind^o$ of projections is uniquely defined.
\end{prop}

\begin{proof} 
By the previous proposition if $D$ has compact resolvents, then $\ind^o(D,A_0,A_1)$ does not depend on the choice of the relative index. Now the proof of Prop. \ref{unique} carries over if we choose the operator $K$ even.
\end{proof}

\subsection{Odd spectral flow and Bott periodicity}
\label{spBottodd}

\begin{ddd}
Let $(D_t)_{t \in [0,1]} \in RSF^-(C([0,1],\hat H_{\A})$ and assume that $D_0$ and $D_1$ are invertible. Then $(C_0((0,1),\hat H_{\A}),\chi(D_t)_{t\in (0,1)})$ is a truly unbounded even Kasparov $(\bbbc,C_0((0,1),\A))$-module. We set
$$\spfl^o((D_t)_{t \in [0,1]}):=\beta[(D_t)_{t \in (0,1)}] \in  K_1(\A) \ .$$ 
\end{ddd}

The definition works also in the case where $\A$ is $\sigma$-unital.

The following properties are easy to verify:

(1) $\spfl^o$ is additive with respect to concatenation of paths.

(2) The restriction of $\spfl^o$ to $BSF^-(C([0,1],\hat H_{\A}))$ is functorial in $\A$.

(3) $\spfl^o((D_t)_{t \in [0,1]})= \spfl^o((f(D_t))_{t \in [0,1]})$ for any odd non-decreasing function $f \in \C(\bbbr)$ with $f(0)=0$ and $f'(0)>0$.

(4) The odd spectral flow is homotopy invariant in a sense analogous to \S \ref{spflone}.

Let $\chi_i \in \C(\bbbr),~i=0,1,$ be such that $\chi_0|_{(-\infty,0])}=1$ and $\chi_0|_{[\frac 13,\infty)}=0$, furthermore $\chi_1|_{[1,\infty)}=1$ and $\chi_1|_{(-\infty,\frac 23]}=0$.

\begin{lem}
Assume that $D \in RSF^-(\hat H_{\A})$ and that $A_0$ and $A_1$ are trivializing operators for $D$. Let $P=1_{\ge 0}(D+A_0)$ and $Q=1_{\ge 0}(D+A_1)$. Then
$$\spfl^o(D + \chi_0(t) A_0+ \chi_1(t) A_1)=\ind^o(Q,P) \in K_1(\A) \ .$$ 
\end{lem} 
 
\begin{proof}
Let $\chi$ be a normalizing function for $D$ such that $\chi(D+A_0)=2P-1$ and $\chi(D+A_1)=2Q-1$. Then $\chi(D + \chi_0(t) A_0+ \chi_1(t) A_1)-((1-t)(2P-1)+ t(2Q-1)) \in C_0((0,1),K(H_{\A}))$. It follows that the class of $[F(\chi(D + \chi_0(t) A_0+ \chi_1(t) A_1))] \in KK_0(\bbbc,\A)$ agrees with the class of $[F((1-t)(2P-1)+t(2Q-1))]$ in $KK_0(\bbbc,C_0((0,1),\A))$. 
\end{proof}

\begin{theorem}
Let $(D_t)_{t\in [0,1]} \in RSF^-(C([0,1],\hat H_{\A}))$ such that $D_0$ and $D_1$ are invertible. Assume that locally trivializing families exist for $(D_t)_{t \in [0,1]}$. Then
$$\spfl_i^o((D_t)_{t\in [0,1]},0,0)= \spfl^o((D_t)_{t\in [0,1]}) \ .$$
\end{theorem}

\begin{proof} The proof uses Cor \ref{uniqueodd} and proceeds as its ungraded counterpart.
\end{proof}

\subsection{Uniqueness of the odd spectral flow}

\begin{ddd}
An odd spectral flow is a map $$\Spfl^o:\{D \in RSF^-(C([0,1],\hat H_{\A}))~|~ D(0) \mbox{ and } D(1) \mbox{ have a bounded inverse} \} \to K_1(\A)$$ with the following properties:
\begin{itemize}
\item[(I)] The restriction of $\Spfl^o$ to $BSF^-(C([0,1],\hat H_{\A}))$ is functorial in $\A$.
\item[(II)] $\Spfl^o(D)= \Spfl^o(\chi(D))$ for any  normalizing function $\chi$ of $D$.
\item[(III)] It is additive with respect to direct sums of operators.
\item[(IV)] The odd spectral flow of a path that is symmetric with respect to the point $\frac 12$ vanishes.
\item[(V)] Let $V$ be a $\bbbz/2$-graded free finitely generated submodule of $\hat H_{\A}$ and let $F_1, F_2$ be odd involutions on $V$. Let $F^{\perp}$ be an odd involution on the orthogonal complement of $V$. Then $\Spfl^o(F^{\perp} \oplus (tF_2+(1-t)F_1))=[U(\frac 12(F_2+1),\frac 12(F_1+1))] \in K_1(\A)$.
\end{itemize}               
\end{ddd}

\begin{prop}
\label{axspecodd}
The odd spectral flow $\Spfl^o$ is uniquely defined and agrees with $\spfl^o$.
\end{prop}

\begin{proof} The proof is very similar to its ungraded counterpart. 
 
(1) The odd spectral flow of a loop $T \in BSF^-(C(S^1,\hat H_{\A}))$ such that $T(0)$ is an involution and $T^2-1$ is compact, only depends on its class $[T] \in KK_1(\bbbc,C_0((0,1)) \ten \A) \subset KK_1(\bbbc,C(S^1) \ten \A)$: 
Assume $[T]=0$. Choose an involution $F \in B^-(\hat H_{\A})$.
Since $[T]=0$, there is an operator $\tilde T \in BSF(C(S^1 \times [0,1],\hat H_{\A}\oplus \hat H_{\A})$ such that $\tilde T(\cdot,0)=T \oplus F$ and $\tilde T(\cdot,1)=T(0) \oplus F$. By functoriality the spectral flow of $T \oplus F$ equals the spectral flow of the loop $T'$ with $T'(x)= \tilde T(0,2x)$ for $x \in [0,\frac 12]$ and $T'(x)=\tilde T(0,1-2x)$ for $x \in [\frac 12,1]$. The spectral flow of $T'$ vanishes by (IV).   

(2) The spectral flow is determined by its action on loops as in (1):

Let $\chi$ be a normalizing function for $D$ such that $\chi(D(0))$ and $\chi(D(1))$ are involutions. Then $P_0:= 2 \chi(D(0))-1$ and $P_1:=2\chi(D(1))-1$ are Lagrangian projections. There is an even unitary $U$ with $UP_0U^*=P_1$. By the contractibility of $\U(H_{\A})$ there is a continuous path of even unitaries $\tilde U:[0,1] \to \U(\hat H_{\A})$ unique up to homotopy with $\tilde U(0)=U,~\tilde U(1)=1$. Define the loop $L(D)$ in $BSF^-(C([0,1],\hat H_{\A}))$ by $L(D)(x)=\tilde U(2x) \chi(D_0) \tilde U(2x)^*$ for $x \in [0,\frac 12]$ and $L(D)(x)= \chi(D(2x-1))$ for $x \in [\frac 12,1]$. By construction $L(0)=L(1)$ is invertible, $L(D)^2-1 \in K(C([0,1],\hat H_{\A}))$ and $\Spfl^o(L(D))= \Spfl^o(D)$.

(3) We evaluate condition (V).  On the loops of the form $F_V \oplus (tF_2+(1-t)F_1)$ the map $\Spfl^o$ agrees with Bott periodicity. Furthermore since Bott periodicity is an isomorphism the classes of these loops generate $KK_0(\bbbc,C_0((0,1),\A))$. Hence  $\Spfl^o(L(D))=\beta[L(D)]$. 

Since $\spfl^o$ fulfills the axioms, $\Spfl^o=\spfl^o$. 
\end{proof}

\subsection{Odd spectral flow and spectral flow}

Here we allow $\A$ to be $\sigma$-unital.

If $D$ is a regular Fredholm operator on $H_{\A}$, then we define $$G(D):=\left(\begin{array}{cc} 0 & D^* \\ D & 0 \end{array}\right) \in RSF^-(\hat H_{\A}) \ .$$

Let $\beta:K_{i+1}(C_0((0,1),\A)) \cong K_i(\A)$ for $i \in \bbbz/2$.

\begin{prop} 
\label{spspodd}
Consider $\cos(\pi x)$ as a multiplication operator on the $C_0((0,1))$-module $C_0((0,1))$.
\begin{enumerate}
\item Let $(D_t)_{t \in [0,1]} \in RSF(C([0,1],H_{\A})$ with $D_0,D_1$ invertible. Then for $(G(D_t - i\cos(\pi x)))_{t \in [0,1]} \in RSF(C([0,1],\hat H_{C_0((0,1),\A)}))$
$$\spfl((D_t)_{t \in [0,1]})=\beta \spfl^o(G(D_t - i\cos(\pi x))_{t \in [0,1]}) \ .$$
\item Let $(D_t)_{t \in [0,1]} \in RSF^-(C([0,1],\hat H_{\A}))$. Then for $(D_t+ \sigma \cos(\pi x))_{t \in [0,1]} \in RSF(C([0,1],\hat H_{C_0((0,1),\A)})$
$$\spfl^o((D_t)_{t \in [0,1]})= \beta\spfl((D_t+ \sigma \cos(\pi x))_{t \in [0,1]}) \ .$$
\end{enumerate}
\end{prop}

\begin{proof}
(1) We use a result which will be proven in the following section. Let $\chi$ be a normalizing function for $(D_t)_{t\in [0,1]}$ . Then Lemma \ref{susp} and Prop. \ref{spbott} imply that the spectral flow of $(D_t)_{t\in [0,1]}$ equals the image of $$[(C_0((0,1)^2,H_{\A}),G(\sin(\pi x) \chi(D_t) -i \cos(\pi x))] \in KK_0(\bbbc, C_0((0,1)^2,\A))$$ under Bott periodicity. Using the fact that for a selfadjoint operator $F$ the operator $F- i\cos(\pi t)$ is invertible for $t\neq \frac 12$ it is easy to check that the right hand side is invariant under the homotopy $s\mapsto G((1-s)\sin(\pi t) \chi(D_t)+ sD_t -i \cos(\pi t))$ of elements in $RSF(\hat H_{C_0((0,1)^2,\A)})$. The results in \S \ref{spBottodd} imply that $\beta\circ \spfl^o$ equals the Bott periodicity map from $KK_0(\bbbc, C_0((0,1)^2,\A))$ to $K_0(\A)$.

(2) is similar and left to the reader.
\end{proof}

{\bf Remark:}  The divisor flow introduced in \cite{m} and its multiparametric odd generalization \cite{lmp} are related to the odd spectral flow as follows: Let $\A=C_0(\bbbr^{2k-1})$, let $M$ be a closed Riemannian manifold and $E$ hermitian vector bundles on $M$. Let $(D_t)_{t \in [0,1]}$ be a continuous path of elliptic elements in $CL^m_{2k-1}(M,E)$, in the notation of \cite{lmp}, and assume that $D_0,D_1$ are invertible. Then $(D_t)_{t\in [0,1]}$ defines a path of regular Fredholm operators on the Hilbert $C_0(\bbbr^{2k-1})$-module $C_0(\bbbr^{2k-1},L^2(M,E))$. The $K$-theoretic description of the divisor flow in \cite[Prop. 2.11]{lmp} implies that
$$\DF((D_t)_{t \in [0,1]})=\spfl^o(G(D_t)_{t \in [0,1]}) \in K_1(C_0(\bbbr^{2k-1})) \cong \bbbz \ .$$  
In a similar way the even divisor flow defined in \cite{lmp} is related to the spectral flow.

The previous proposition reproves and generalizes \cite[Prop. 3.1]{lmp}. The original proof uses $\eta$-invariants and cannot be generalized to a $C^*$-algebraic context.

\section{A splitting formula}
\label{split}

As an application we explain the role of the noncommutative Maslov index in a splitting formula for the spectral flow. The splitting formula for the spectral flow is derived from a splitting formula for the index for families of Dirac operators. We only discuss the case of vertical Dirac operators on fiber bundles with odd-dimensional fibers but the proofs carry over to the case of even-dimensional fibers and to Dirac operators over $C^*$-algebras. 

Consider a fiber bundle $\pi:M \to B$ with closed odd-dimensional fibers. Endow it with a vertical Riemannian metric. Let 
$E$ be a vector bundle on $M$ with a vertical Dirac bundle structure and let $D=(D_b)_{b \in B}$ be the vertical Dirac operator associated to $E$. 

We define the Hilbert $C(B)$-module $L^2_v(M,E)$  as the completion of $C(M,E)$ with respect to the $C(B)$-valued scalar product induced by the $L^2$-scalar product on the fibres. Note that $L^2_v(M,E)$ is isomorphic to $H_{C(B)}$. We denote the closure of $D$ on $L^2_v(M,E)$ with domain $\C(M,E)$ by $D$ as well.

Let $N \subset M$ be a fiber bundle of hypersurfaces such that $M=M_0 \cup_N M_1$ for two manifolds with boundary $M_0,M_1 \subset M$ with $\ra M_i=N$. We assume that all vertical structures are of product type near the boundary. Let $x$ be the normal coordinate in a neighborhood of $N$ with $x(N)=0$ and $x|_{M_0} \le 0$. Denote $E_N=E|_N$. The operator $ic(dx):E_N \to E_N$ defines a grading on $E_N$. Near $N$ the Dirac operator takes the form 
$c(dx)(\ra_x + D_N)$ where $D_N$ is a vertical Dirac operator associated to $E_N$ on $N$. We have that $c(dx)D_N=-D_N c(dx)$. 

Let $x_0$ resp. $x_1$ be the boundary defining coordinates near the boundary of $M_0$ resp. $M_1$, hence $x_0=x$ on $M_0$ and $x_1=-x$ on $M_1$. We denote the restriction of $D$ to $M_0$ resp. $M_1$ by $D_0$ resp. $D_1$. 

We define generalized Atiyah-Patodi-Singer boundary conditions for $D_0$ and $D_1$. Near the boundary of $M_0$ resp. $M_1$ the operator $D$ equals 
$c(dx_0)(\ra_{x_0} + D_N)$ resp. $c(dx_1)(\ra_{x_1} - D_N)$. We choose smoothing odd trivializing operators $A_0$ resp. $A_1$ of $D_N$ and write $P_i=1_{\ge 0}(D_N+A_i)$. Define the fiberwise domains $$\dom D_{0b}:= \{f \in \C(M_{0b},E|_{M_0})~|~ P_{0b}(f|_{N_b})=0\}$$ and 
$$\dom D_{2b}:= \{f \in \C(M_{2b},E|_{M_1})~|~ (1-P_{1b})(f|_{N_b})=0\} \ .$$ 
Define now $\dom D_0$ resp. $\dom D_1$ as in \S \ref{fam}. 

We consider $N\times [0,1]$ as a fiber bundle with base space $B$ and denote the pullback of $E_N$ to $N\times [0,1]$ by $E_N$ as well. Let $D_Z$ be the closure of the operator $c(dx)(\ra_x + D_N)$ on $L^2_v(N\times [0,1],E_N)$ with domain
$$\{f \in C(N\times [0,1],E_N))~|~ f \mbox{ is fiberwise smooth and } \dots $$ 
$$ \dots (1-P_0)(f|_{N\times\{0\}})=P_1(f|_{N\times \{1\}})=0\} \ .$$

The operators $D,D_0,D_1,D_Z$ are fiberwise selfadjoint regular and Fredholm. Furthermore they have compact resolvents. It follows from Lemma \ref{critreg2} as in the proof of Prop. \ref{boundval} that $D,D_0,D_1,D_Z$ are regular, selfadjoint and Fredholm on the respective Hilbert $C(B)$-modules. 

Hence $D,D_0,D_1,D_Z$ define truly unbounded odd Kasparov $(\bbbc,C(B))$-modules. 
We have that 
$$\ind^o(D_0) + \ind^o(D_1) = \ind^o(D) + \ind^o(D_Z) \in K^1(B) \ .$$
The proof of this formula follows from a simple modification of the $K$-theoretic relative index theorem \cite{bu} in order to include the boundary conditions.

\begin{prop}
$$\ind^o(D_Z)=\ind^o(P_0,P_1) \ .$$ 
\end{prop}

\begin{proof}
Let ${\cal D}=i\sigma(\ra_t + D_N +\chi_0(t) A_0 + \chi_1(t) A_1)$.

In view of Prop. \ref{indspecsecodd} we have to show that
$\ind^o(D_Z) = \ind^o({\cal D}).$

The proof is a relative $K$-theoretic index theorem. Note first that both sides vanish for $A_0=A_1$. We may assume that $P_0,P_1$ are spectral sections. Since $D_N$ admits spectral sections there is a spectral section $R$ with $RP_1=R$ and $RP_2=R$. Then $P_1-R$ and $P_2-R$ are projections onto Lagrangian subspaces in 
$\Ran(1-\sigma R \sigma-R)$. Let $D_N'=D_N-(1-\sigma R \sigma-R)D_N(1-\sigma R \sigma-R)$. We denote by $D_Z'$ resp. ${\cal D}'$ the operators obtained from $D_Z$ resp. ${\cal D}$ by replacing $D_N$ with $D_N'$. Clearly $\ind^o(D_Z)=\ind^o(D_Z')$ and $\ind^o({\cal D})=\ind^o({\cal D}')$. Let $U$ be the even unitary that equals the identity on $\Ran(\sigma R \sigma+R)$ and equals $U(R-P_0,R-P_1)$ on $\Ran(1-\sigma R \sigma-R)$. Then $UP_0U^*=P_1$ and $UD_N'U^*=D_N'$. Now the proof proceeds as the proof of Cor. \ref{standiso2}.  
\end{proof}

It follows that
$$\ind^o(D_0) + \ind^o(D_1) = \ind^o(D) + \ind^o(P_0,P_1) \in K^1(B) \ .$$

A splitting formula for the spectral flow for a loop of families of Dirac operators follows immediately: We consider the loop as a family of Dirac operators with base space $S^1 \times B$. Let the notation be as before, now with base space $S^1 \times B$, and assume that the restriction of $D$, $D_0$, $D_1$ to the fiber of $1 \in S^1$ has vanishing index in $K^1(B)$. Then the previous formula holds   
in $K^1((0,1) \times B) \subset K^1(S^1 \times B)$. By applying the Bott periodicity map $\beta:K^1((0,1) \times B) \to K^0(B)$ to both sides we get
$$\spfl(D_0) + \spfl(D_1) = \spfl(D) + \beta (\ind^o(P_0,P_1)) \in K^0(B) \ .$$ 

The definition of the odd relative index of projections $\ind^o_1$ and Cor. \ref{maslbott} justify to call $\beta \ind^o(P_0,P_1)$ the Maslov index of $(P_0,P_1)$. If the kernel of $D_N$ is a trivial vector bundle over $S^1 \times B$, then we may choose $P_i= 1_{>0}(D_N) + P_{L_i}$, where $P_{L_i}:\Ker D_N \to L_i$ is a projection onto a Lagrangian submodule $L_i \subset \Ker D_N$. We assume that $L_0,L_1$ are transverse at $1 \in S^1$. Then  $$\beta \ind^o(P_0,P_1)=\mu(L_0,L_1) \in K_0(B) \ .$$

In \cite{n} a splitting formula for the spectral flow for general paths of Dirac type operators was proven. In contrast to the proof in \cite{n}, which uses the unique continuation property of Dirac type operators, our proof generalizes to elliptic operators which are not of Dirac type and to fiber bundles with noncompact fibers under the condition that $N$ is compact and the operators have the above form near $N$. 

The splitting index formulas \cite[Prop. 11]{mp2} (for families) and \cite[Theorem 9]{lp2} (for Dirac operators over $C^*$-algebras) use a different definition for the odd index of a selfadjoint operator. We will explain this alternative definition in the following and prove that it agrees with our definition. By comparing the splitting theorems one obtains that the difference elements of Lagrangian spectral sections defined in \cite{mp} and \cite{lp2} and our definition are equivalent.   

Let $\A$ be a $\sigma$-unital $C^*$-algebra and let $F \in B(H_{\A})$ be selfadjoint with $F^2-1 \in K(H_{\A})$. The operator $$F_t:=\left(\begin{array}{cc} 0 & \cos(\pi t)- i\sin(\pi t)F \\ \cos(\pi t)+ i\sin(\pi t)F & 0 \end{array}\right)$$ on the $\bbbz/2$-graded Hilbert $C_0((0,1),\A)$-module $C_0((0,1),\hat H_{\A})$ defines an element $[F_t]$ in $KK_0(\bbbc,C_0((0,1),\A)) \cong K_0(C_0((0,1),\A)$. Via Bott periodicity we obtain an element $\Ind^o(F) \in K_1(\A)$. This construction is a generalization of a suspension map introduced in \cite{as}. The following lemma implies that $\Ind^o(F)$ agrees with the image of $[F] \in KK_1(\bbbc,\A)$ under the standard isomorphism $KK_1(\bbbc,\A) \cong K_1(\A)$. In order to define an index $\Ind^o(D)$ for $D\in RSF(H_{\A})$ we choose a normalizing function $\chi$ and apply this construction to $\chi(D)$. Then the following lemma also implies that $\Ind^o(D):=\Ind^o(\chi(D))$ agrees with $[D] \in KK_1(\bbbc,\A)$ under the standard isomorphism $KK_1(\bbbc,\A)\cong K_1(\A)$. In particular $\Ind^o(F)=\ind^o(F)$ and $\Ind^o(D)=\ind^o(D)$. 

\begin{lem} 
\label{susp}
The class $[(F_t)_{t\in [0,1]}] \in KK_0(\bbbc,C_0((0,1),\A))$ corresponds to $[F] \in KK_1(\bbbc,\A)$ under Bott periodicity.
\end{lem}

\begin{proof}
The Bott element in $KK_1(\bbbc,C_0((0,1)))$ is represented by the operator $-\cos(\pi t)$ on the $C_0((0,1))$-module $C_0((0,1))$. In the following we show that 
$$[- \cos(\pi t)] \ten [F] = [F_t] \ .$$ Let $C_1^1,C_1^2$ be copies of the Clifford algebra $C_1$ and let $c_i \in C_1^i$ be odd, selfadjoint with $c_i^2=1$. Recall that the Kasparov product
$$KK_1(\bbbc,C_0((0,1)) \ten KK_1(\bbbc,\A) \to KK_0(\bbbc,C_0((0,1),\A))$$ is defined via the Kasparov product
$$KK_0(\bbbc, C_1^1\ten C_0((0,1))) \ten KK_0(\bbbc, C_1^2\ten \A) \to KK_0(\bbbc,C_1^1\ten C_1^2 \ten C_0((0,1),\A)) \ .$$
The computation is now straight-forward:
Under the isomorphism $KK_1(\bbbc,C_0((0,1)))\cong KK_0(\bbbc, C_1^1\ten C_0((0,1)))$ resp. $KK_1(\bbbc,\A)\cong KK_0(\bbbc, C_1^2\ten \A)$ the element $[-\cos(\pi t)]$ resp. $[F]$ corresponds to  $[C_1^1 \ten C_0((0,1)),-c_1\cos(\pi t)]$ resp. $[C_1^2 \ten H_{\A}, c_2 F]$. We calculate the Kasparov product of the last two elements in $KK_0(\bbbc,C_1^1\ten C_1^2 \ten C_0((0,1),\A))$. By tensoring $c_2F$ with the identity we obtain an operator on the Hilbert $C_1^1\ten C_1^2 \ten C_0((0,1),\A)$-module $C_1^1\ten C_1^2 \ten C_0((0,1),H_{\A})$ which we denote by $c_2F$ as well. It is clear that $c_2F$ is a $c_2F$-connection for $C_1^1 \ten C_0((0,1))$ (see \cite[Def. 18.3.1]{bl} ). Then the operator $-c_1\cos(\pi t)+ c_2\sin(\pi t)F$ on  $C_1^1\ten C_1^2 \ten C_0((0,1),H_{\A})$ represents the Kasparov product  $[-c_1\cos(\pi t)]] \ten [c_2 F] \in KK_0(\bbbc,C_1^1\ten C_1^2 \ten C_0((0,1),\A))$ (see \cite[Ex. 18.4.2.d]{bl}). Now $p= \frac 12(1 + ic_1c_2) \subset C_1^1 \ten C_1^2$ is an even projection of rank one, hence the homomorphism $f:\bbbc \to C_1^1 \ten C_1^2,~x\mapsto xp$ induces the standard isomorphism in $KK$-theory. Let $M_2(\bbbc)$ be endowed with the grading induced by $M_2(\bbbc)\cong \End(\bbbc^+ \oplus \bbbc^-)$ and identify $$C_1^1 \ten C_1^2 \cong M_2(\bbbc),~c_1 \mapsto \left(\begin{array}{cc} 0 &-1 \\ -1 & 0 \end{array}\right) \mbox{ and } c_2 \mapsto \left(\begin{array}{cc} 0 & -i \\ i & 0 \end{array}\right) \ .$$ 
Then $p=\left(\begin{array}{cc} 1 & 0\\0 &0 \end{array}\right)$. One verifies that on the Hilbert $C_0((0,1),\A)$-module $(C_1^1\ten C_1^2 \ten C_0((0,1),H_{\A})) \ten_f (\bbbc^+ \oplus \bbbc^-)$ the operator $F_t$ equals $(-c_1\cos(\pi t)+ \sin(\pi t)c_2F) \ten_f 1$, hence $[F_t]\in KK_0(\bbbc,C_0((0,1),\A))$ is the image of the Kasparov product $[-c_1\cos(\pi t)]] \ten [c_2 F] \in KK_0(\bbbc,C_1^1\ten C_1^2 \ten C_0((0,1),\A))$ under the standard isomorphism.
\end{proof}

\textsc{Nieders\"achsische Landesbibliothek\\
Gottfried Wilhelm Leibniz Bibliothek\\
Waterloostr. 8\\
30169 Hannover\\
Germany} 

\textsc{Email: ac.wahl@web.de}

\end{document}